\documentclass{article}
\usepackage[american]{babel}
\usepackage{amsfonts,amsmath,amssymb}
\usepackage{xypic}
\xyoption{all}
\input{xypic}
\newdir{ >}{{}*!/-8pt/\dir{>}}

\input liemacs.sty
\def\P{{\mathbb P}}

\def\k{\mathfrak k}

\begin{document}
\title{Holomorphic current groups -- Structure and Orbits} 

\author{Martin Laubinger\\
        Allianz Deutschland AG\\
        martin.laubinger@gmail.com\\
\and    Friedrich Wagemann \\
        Universit\'e de Nantes\\
        wagemann@math.univ-nantes.fr}

\maketitle


\begin{abstract}
Let $K$ be a finite-dimensional, 1-connected complex Lie group,
and let $\Sigma_k=\Sigma\setminus\{p_1,\ldots,p_k\}$ be a compact 
connected Riemann surface $\Sigma$, from which we have extracted $k\geq 1$ 
distinct points.   
We study in this article the regular Fr\'echet-Lie group 
${\mathcal O}(\Sigma_k,K)$ of holomorphic maps from $\Sigma_k$ to $K$
and its central extension $\widehat{{\mathcal O}(\Sigma_k,K)}$. 
We feature especially the automorphism groups of these Lie groups
as well as the coadjoint orbits of $\widehat{{\mathcal O}(\Sigma_k,K)}$
which we link to flat $K$-bundles on $\Sigma_k$. 
\end{abstract}

\section*{Introduction}

Let $K$ be a connected Lie group. 
The {\it loop group} $LK:={\mathcal C}^{\infty}(\bS^1,K)$ is an important 
symmetry group 
in physics, and its mathematical study has led to considerable research
with respect to its structure, its representations and their applications. 
As physicists usually prefer the symmetry Lie algebra to the group, the 
loop Lie algebra $L{\mathfrak k}:={\mathcal C}^{\infty}(\bS^1,{\mathfrak k})$ 
where ${\mathfrak k}$ is the Lie algebra of $K$, also plays an important role. 
Most of the algebraic theory is based on $L{\mathfrak k}$ and its central 
extension $\widehat{L{\mathfrak k}}$ which is closely related to affine 
Kac-Moody Lie algebras. 
There is interesting 
geometry coming about when lifting the central extension 
$\widehat{L{\mathfrak k}}$ to group level defining the $\widehat{LK}$. 
A standard reference about the loop group is Pressley-Segal's book 
\cite{PreSeg}. 

Searching for similar mapping groups in higher dimensions, Etingof and Frenkel
\cite{EtiFre} came up with the group $K^{\Sigma}:={\mathcal C}^{\infty}(\Sigma,K)$
 for a 
compact connected Riemann surface $\Sigma$, and studied its central extensions,
its automorphisms and coadjoint orbits. It turns out that the universal central 
extension has infinite dimensional center, but there is an interesting 
sub-extension by the Jacobian of $\Sigma$. The orbits of the coadjoint 
representation for 1-connected $K$ correspond bijectively to equivalence 
classes of holomorphic principal $K$-bundles over the surface $\Sigma$.
These bundles carry by construction a flat connection.   

On the other hand, more classes of examples arise from mapping spaces with
different regularity conditions. One important class of examples are 
spaces of meromorphic maps. Here the research mainly restricted to the Lie 
algebra side (it is difficult or even impossible to define Lie groups for 
these Lie algebras), namely to {\it Krichever-Novikov algebras}. 
Given a compact connected Riemann surface $\Sigma$ and $k\geq 1$ distinct points
$p_1,\ldots,p_k$, the open Riemann surface 
$$\Sigma_k\,:=\,\Sigma\setminus\{p_1,\ldots,p_k\}$$
is in fact an affine algebraic curve. 
Krichever-Novikov \cite{KN1}, \cite{KN2}, \cite{KN3} for the genus zero case 
and Schlichenmaier, Schlichenmaier-Sheinman (see \cite{Sch}) for arbitrary 
genus define meromorphic 
analogues of the loop algebras as $Reg(\Sigma_k)\otimes_{\C}{\mathfrak k}$
for any complex Lie algebra ${\mathfrak k}$, where $Reg(\Sigma_k)$ is the 
algebra of regular functions on the algebraic variety $\Sigma_k$. They also 
study central extensions, automorphisms, coadjoint orbits and representations 
of these Lie algebras. 

The present work takes its place between Etingof-Frenkel's work and the work on
Krichever-Novikov algebras. Our goal is to study central extensions, 
automorphisms and coadjoint orbits of holomorphic analogues of 
Krichever-Novikov algebras. A new feature is that we have Fr\'echet-Lie groups 
to these Lie algebras here. 
Namely, instead of considering only {\it regular} 
functions on the affine algebraic curve $\Sigma_k$, one may consider all 
{\it holomorphic maps} on the Stein manifold $\Sigma_k$. These latter may have 
essential singularities in the points $p_1,\ldots,p_k$, while the former are
meromorphic, i.e. have at most poles of finite order in these points. 
We are thus considering here the topological group ${\mathcal O}(\Sigma_k,K)$
of holomorphic maps form $\Sigma_k$ to $K$ in the topology of uniform 
convergence on compact sets.  
In fact, ${\mathcal O}(\Sigma_k,K)$
is an infinite-dimensional Fr\'echet-Lie group \cite{NeeWag} with Lie algebra
${\mathcal O}(\Sigma_k,{\mathfrak k})$, the Lie algebra of holomorphic maps 
from $\Sigma_k$ to the Lie algebra ${\mathfrak k}$. 

In our study of the structure of ${\mathcal O}(\Sigma_k,K)$ and its central 
extension $\widehat{{\mathcal O}(\Sigma_k,K)}$ which occupies Section 2,
we build on our previous 
work concerning the holomorphic 
current algebra \cite{NeeWag}, \cite{NeeWag3} and on previous studies of
of the infinite-dimensional Lie theory of current groups \cite{Neeb},
\cite{NeeWag1}, \cite{NeeMai} part of which we adapt to the holomorphic 
setting. A new result is the lifting of the Lie group structure to
the central extension (Corollary 2.13). 

In Section 3, we describe the automorphism groups of 
${\mathcal O}(\Sigma_k,K)$ and of its central extension.
Concrete results are the computation of the 
automorphism group (Corollary \ref{Aut(G)}). We start the problem 
of determining the automorphism groups of the central extensions 
(Proposition 3.9), but leave it to further study to determine these groups 
explicitly for the different complex simple Lie algebras ${\mathfrak k}$.
  
Section 4 is the heart of the present article. Here we study the coadjoint 
orbits of ${\mathcal O}(\Sigma_k,K)$ in the smooth dual of 
$\widehat{{\mathcal O}(\Sigma_k,{\mathfrak k})}$ and establish relations to 
flat principal $K$-bundels on $\Sigma_k$ as well as to the coadjoint orbits 
of the loop group $LK$ (Propositions 4.3 and 4.8). Section 1 prepares the 
necessary material about flat bundles. 
As usual, these coadjoint orbits carry the 
{\it Kostant-Kirillov-Souriau} symplectic form, and they are thus examples
of (weakly) holomorphic symplectic manifolds (see Corollary 4.12).

\vspace{.5cm}

\noindent{\bf Acknowledgements:}
Both authors thank Karl-Hermann Neeb for constant support, very useful 
discussions and references. ML thanks Laboratoire de Math\'ematiques 
Jean Leray which financed his stay in Nantes in 2010 where part of 
this work was done.

\section{Holomorphic $1$-forms and flat principal bundles}

\subsection{Holomorphic $1$-forms and the logarithmic derivative}
\label{section_flat_bundles}

In this subsection, we discuss the conditions under which a Lie algebra valued
holomorphic $1$-form is the logarithmic derivative of a group-valued 
holomorphic function. We follow closely \cite{NeeWag}. 

Let $K$ be a (possibly infinite dimensional) regular complex Lie group.  
For a complex manifold $M$, let us denote by ${\mathcal O}(M,K)$ the group of
holomorphic maps from $M$ to $K$ and by $\Omega^1(M,{\mathfrak k})$ the 
space of holomorphic $1$-forms on $M$ with values in the Lie algebra
${\mathfrak k}$. Let $\kappa$ be the Maurer-Cartan form on $K$, 
i.e. the 
unique holomorphic left invariant $1$-form with values in ${\mathfrak k}$ on 
$K$ corresponding to $\id_{\mathfrak k}$ in 
$\Hom({\mathfrak k},{\mathfrak k})$.
For $f\in{\mathcal O}(M,K)$, we denote by 
$$\delta(f)\,:=\,f^*(\kappa)\,=:\,f^{-1}df\,\in\Omega^1(M,{\mathfrak k})$$
the {\it (left) logarithmic derivative of $f$}. We obtain a map
$$\delta:{\mathcal O}(M,K)\to \Omega^1(M,{\mathfrak k})$$
satisfying the cocycle condition
$$\delta(f_1f_2)\,=\,{\rm Ad}(f_2)^{-1}\delta(f_1)+\delta(f_1).$$
From this it follows easily that if $M$ is connected, then
$$\delta(f_1)=\delta(f_2)\,\,\,\,\,\Leftrightarrow\,\,\,\,\,(\exists k\in K)
\,\,\,f_2=\lambda_k\circ f_1,$$
where $\lambda_k$ denotes left multiplication by $k\in K$. If $K$ is abelian,
then $\delta$ is a group homomorphism whose kernel consists of the locally 
constant maps $M\to K$. There is a (right) action of ${\mathcal O}(M,K)$ on 
$\Omega^1(M,{\mathfrak k})$ given by
\begin{equation}   \label{gauge_action}
\alpha * f \,:=\,\delta(f) + {\rm Ad}(f)^{-1}\alpha,
\end{equation}
derived from the above cocycle condition. We will call this action the {\it 
gauge action}.  

We call $\alpha\in\Omega^1(M,{\mathfrak k})$ {\it integrable} if there exists
a holomorphic function $f:M\to K$ with $\delta(f)=\alpha$. We say that $\alpha$
is {\it locally integrable} if each point $m\in M$ has an open neighborhood
$U$ such that $\alpha|_U$ is integrable. In order to describe conditions for
the integrability of an element $\alpha\in\Omega^1(M,{\mathfrak k})$, we 
define the bracket 
$$[.,.]:\Omega^1(M,{\mathfrak k})\times\Omega^1(M,{\mathfrak k})\to
\Omega^2(M,{\mathfrak k})$$
by
$$[\alpha,\beta]_p(v,w)\,:=\,[\alpha_p(v),\beta_p(w)]-[\alpha_p(w),\beta_p(v)]
\,\,\,\,\,{\rm for}\,\,\,\,\,v,w\in T_p(M).$$ 
Note that $[\alpha,\beta]=[\beta,\alpha]$. By definition, the {\it 
Maurer-Cartan space} associated to the complex manifold $M$ and the complex Lie
algebra ${\mathfrak k}$ is the space of closed non-abelian differential 
$1$-forms
$$Z^1_{\rm dR}(M,{\mathfrak k})\,:=\,\{\alpha\in \Omega^1(M,{\mathfrak k})\,|\,
d\alpha+\frac{1}{2}[\alpha,\alpha]=0\,\}.$$

The following theorem can be found with a full proof in \cite{NeeWag}. 
Here we denote by
$${\rm evol}_K:{\mathcal C}^{\infty}([0,1],{\mathfrak k})\to K$$
the {\it evolution map} which associates to the initial value problem 
$$\gamma'(t)\,=\,\gamma(t)\cdot\xi$$
associated to
$\xi\in{\mathcal C}^{\infty}([0,1],{\mathfrak k})$ its solution 
$\gamma_{\xi}(1)\in K$ at $1\in K$. 
The Lie group $K$ is called {\it regular} if for all 
$\xi\in{\mathcal C}^{\infty}([0,1],{\mathfrak k})$, the solution 
$\gamma_{\xi}$ exists and the evolution map is smooth. 
For a complex regular Lie group $K$, the evolution map is
holomorphic, cf {\it loc. cit.} Lemma $3.4$.

\begin{Theorem} \label{fundamental_theorem}
Let $M$ be a complex manifold, $K$ be a regular complex Lie group and 
$\alpha\in\Omega^1(M,{\mathfrak k})$.
\begin{itemize}
\item[(1)] $\alpha$ is locally integrable if and only if 
$\alpha\in Z^1_{\rm dR}(M,{\mathfrak k})$.
\item[(2)] If $M$ is $1$-connected and $\alpha$ is locally integrable, then 
$\alpha$ is integrable.
\item[(3)] Suppose that $M$ is $1$-connected, fix $m_0\in M$ and let 
$\alpha\in Z^1_{\rm dR}(M,{\mathfrak k})$. 
Using piecewise smooth representatives 
of homotopy classes, we obtain a well-defined group homomorphism
$${\rm per}_{\alpha}^{m_0}:\pi_1(M,m_0)\to K,\,\,[\gamma]\mapsto{\rm evol}_K
(\gamma^*\alpha),$$
and $\alpha$ is integrable if and only if this homomorphism is trivial.
\end{itemize}
\end{Theorem}

Sometimes, we will suppress the base point in the notation and write more 
simply ${\rm per}_{\alpha}$ for ${\rm per}_{\alpha}^{m_0}$.  

\subsection{The classification of flat holomorphic principal bundles}

In this subsection, we follow closely \cite{Nee04} and \cite{Lau}. The main
result (Proposition \ref{bundle_exact_sequence}) of this subsection is the 
exact sequence linking solutions of the 
Maurer-Cartan equation, homomorphisms from the fundamental group $\pi_1$ 
to the Lie group $K$ 
and flat principal $K$-bundles, which appears in \cite{Nee04} and \cite{Lau}.
The main idea is to describe flat bundles as 
bundles which become trivial when pulled back to the universal covering. 

Let $M$ be a connected complex manifold and $K$ be a connected complex Lie 
group.
Recall that the universal covering $\widetilde{M}$ of $M$ has a natural 
structure of a complex manifold.   
Let us briefly recall non-abelian $1$-cocycles and non-abelian $1$-cohomology.

Let $\Gamma$ be a group. A {\it $\Gamma$-group} is a group $G$ endowed with 
an action of $\Gamma$ by automorphisms. We denote this action by 
$(\gamma,g)\mapsto\gamma\cdot g$. Given a $\Gamma$-group $G$, a 
{\it $G$-valued $1$-cochain} of $\Gamma$ is simply a map $\Gamma\to G$.
A cochain $c:\Gamma\to G$ is called a {\it non-abelian $1$-cocycle} in case
for all $\alpha,\beta\in\Gamma$, we have
$$c(\alpha)\,\,(\alpha\cdot c(\beta))\,\,c(\alpha\beta)^{-1}\,=\,1.$$
Such a $1$-cocycle is also called a {\it crossed homomorphism}.
On the set $Z^1(\Gamma,G)$ of non-abelian $1$-cocycles, there is a 
$G$-action given by
$$(g\cdot c)(\gamma)\,=\,g\,\,c(\gamma)\,\,(\gamma\cdot g)^{-1}.$$
The {\it non-abelian cohomology set} $H^1(\Gamma,G)$ is by definition the set
of $G$-orbits in $Z^1(\Gamma,G)$. In case $\Gamma$ acts trivially on $G$, 
we have 
$$Z^1(\Gamma,G)\,=\,\Hom(\Gamma,G)\,\,\,\,{\rm and}\,\,\,\,
H^1(\Gamma,G)=\Hom(\Gamma,G)\,/\,G,$$
where $G$ acts on $\Hom(\Gamma,G)$ by the conjugation action in the image:
$(g\cdot c)(\gamma)\,=\,g\,c(\gamma)\,g^{-1}$. 

In the following, we will apply this formalism to the special case where
$\Gamma=\pi_1(M)$ and $G={\mathcal O}(\widetilde{M},K)$, where 
$\widetilde{M}$ is the universal covering of $M$ and $\pi_1(M)$ acts on 
$G$ via its action on $\widetilde{M}$. 

\begin{Proposition}
Let $q_M:\widetilde{M}\to M$ be the universal covering projection. 
To a non-abelian $1$-cocycle $\phi:\pi_1(M)\to {\mathcal O}(\widetilde{M},K)$, 
one associates the holomorphic principal $K$-bundle
$$P_{\phi}\,:=\,(\widetilde{M}\times K)\,/\,\pi_1(M),\,\,\,\,\,[(m,k)]:=
\pi_1(M)\cdot(m,k),$$
where $\pi_1(M)$ acts on the trivial bundle $\widetilde{M}\times K$ by
$$\gamma\cdot(m,k)\,=\,(\gamma\cdot m,\phi(\gamma)(m)k).$$
Then the bundle $q^*_MP_{\phi}$ is holomorphically trivial. 
Conversely, each holomorphic principal $K$-bundle 
$q:P\to M$ for which $q_M^*P$ is holomorphically trivial is equivalent 
to some $P_{\phi}$. For non-abelian $1$-cocycles $\phi,\psi$, we have
$$P_{\phi}\cong P_{\psi}\,\,\,\,\,\Leftrightarrow\,\,\,\,\,(\exists f\in
{\mathcal O}(\widetilde{M},K):\,\,\,\,\psi(\gamma)=f\cdot\phi(\gamma)
\cdot(\gamma\cdot f)^{-1}\,\,\,\,\,{\rm for}\,\,\,\,\,{\rm all}\,\,\,\,\,
\gamma\in\pi_1(M).$$ 
Let ${\rm Bun}_{\widetilde{M}}(M,K)$ denote the set of equivalence classes of
holomorphic principal $K$-bundles whose pullback to $\widetilde{M}$ is 
holomorphically trivial. Then the map
$$H^1(\pi_1(M),{\mathcal O}(\widetilde{M},K))\to
{\rm Bun}_{\widetilde{M}}(M,K),\,\,\,\,\,[\phi]\mapsto[P_{\phi}]$$
is a bijection.
\end{Proposition}

\begin{Proof}
The proof can be adapted from \cite{Nee04} and appears in \cite{Lau}
Lemma 5. Let us recall here only its main features:
The bundle $P_{\phi}$ is holomorphic, because the values of $\phi$ are 
holomorphic functions and $K$ is a complex Lie group. Now
$$q_M^*{P_{\phi}}\,=\,\{([m,k],m)\in P_{\phi}\times\widetilde{M}\,|\,k\in K\},
$$ 
and $m\mapsto([m,1],m)$ is a holomorphic global section. On the other hand, one 
shows that an arbitrary $P$ may be identified with the set 
$q_M^*{P}\,/\,\pi_1(M)$ of $\pi_1(M)$-orbits in $q_M^*{P}$. Furthermore, this
kind of action by biholomorphisms on a trivial principal bundle is always 
given by a non-abelian $1$-cocycle $\phi$, and one has $P\cong P_{\phi}$.

Then, an equivalence $\Phi$ of bundles $P_{\phi}\cong P_{\psi}$ is necessarily 
of the form $\Phi(m,k)=(m,f(m)k)$ for some $f:\widetilde{M}\to K$. The 
compatibility with the $\pi_1(M)$-action is equivalent to the equation 
$\psi(\gamma)=f\cdot\phi(\gamma)\cdot(\gamma\cdot f)^{-1}$. But this 
equation just means that $\phi$ is sent to $\psi$ acting on it by 
$f\in{\mathcal O}(\widetilde{M},K)$, i.e. they belong to the same 
class in $H^1(\pi_1(M),{\mathcal O}(\widetilde{M},K))$.     
\end{Proof} 

\begin{rem}
Suppose for the sake of this remark that $K$ is connected and abelian. 
Note that the map $\phi\mapsto P_{\phi}$ appears in the spectral sequence 
arising in equivariant cohomology for the covering
$$\pi_1(M)\to\widetilde{M}\to M.$$
More precisely, the discrete group $H:=\pi_1(M)$ gives rise to a classifying 
space $BH$ and the total space $EH$ of the universal $H$-bundle. The 
universal covering space $\widetilde{M}$ can be described up to homotopy 
equivalence as $EH\times\widetilde{M}$ and $M$ as 
$(EH\times\widetilde{M})\,/\,H$. The fibration 
$$EH\times\widetilde{M}\to(EH\times\widetilde{M})\,/\,H$$
has typical fiber $EH$ which is contractible, thus we obtain for the 
equivariant cohomology 
$$H^*_H(\widetilde{M},\underline{K})\,=\,H^*((EH\times\widetilde{M})\,/\,H,
\underline{K})\,=\,H^*(M,\underline{K}).$$
On the other hand, the fibration 
$$(EH\times\widetilde{M})\,/\,H\to EH\,/\,H=BH$$ 
gives rise to a spectral sequence in equivariant cohomology whose $E_2$-term is
$H^*(BH,H^*(\widetilde{M},\underline{K}))$ and which converges to
the equivariant cohomology of $\widetilde{M}$: 
$$E_{\infty}^*\,=\,H_H^*(\tilde{M},\underline{K})\,=\,H^*(M,\underline{K}).$$
For non-abelian $K$, these constructions only work for low-dimensional 
(\v{C}ech) cohomology groups.    
We therefore obtain in this framework maps 
$$H^1(B(\pi_1(M)),\check{H}^0(\widetilde{M},\underline{K}))\,=\,
H^1(\pi_1(M),{\mathcal O}(\widetilde{M},K))\stackrel{\partial}{\to}
\check{H}^1(M,\underline{K})\stackrel{q_M^*}{\to}
\check{H}^1(\widetilde{M},\underline{K}).$$
Now translate $\check{H}^1(M,\underline{K})={\rm Bun}(M,K)$
and $\ker(q_M^*)={\rm Bun}_{\widetilde{M}}(M,K)$.
The upshot is that $\partial$ has values in $\ker(q_M^*)$ by exactness, 
and induces a bijection. 
\end{rem}

\begin{Corollary} \label{preceeding_corollary}
The bundle $P_{\phi}$ is trivial if and only if there exists a holomorphic 
function $f\in{\mathcal O}(\widetilde{M},K)$ with  
$$\phi(\gamma)=f\cdot(\gamma\cdot f)^{-1}.$$
For $\phi\in\Hom(\pi_1(M),K)$, this assertion is equivalent to the existence of
an $\alpha\in Z^1_{\rm dR}(M,{\mathfrak k})$ with ${\rm per}_{\alpha}=\phi$.
\end{Corollary}

\begin{Proof}
It remains to show that for $\phi\in\Hom(\pi_1(M),K)$, the existence of 
$f\in{\mathcal O}(\widetilde{M},K)$ with 
$\phi(\gamma)=f\cdot(\gamma\cdot f)^{-1}$ for all $\gamma\in\pi_1(M)$ is
equivalent to the existence of $\alpha\in Z^1_{\rm dR}(M,{\mathfrak k})$ with 
${\rm per}_{\alpha}=\phi$. But the condition 
$\phi(\gamma)=f\cdot(\gamma\cdot f)^{-1}$ may be rephrased as 
$f(\gamma\cdot m)=\phi(\gamma)f(m)$ for all $m\in\widetilde{M}$, 
$\gamma\in\pi_1(M)$. By Theorem \ref{fundamental_theorem}, as the universal 
cover is  $1$-connected, $q_M^*(\alpha)$ is a logarithmic derivative and the 
equation $f(\gamma\cdot m)=\phi(\gamma)f(m)$ means that  
${\rm per}_{\alpha}=\phi$. This argument works both ways.
\end{Proof}

\begin{Proposition}   \label{preceeding_proposition}
A holomorphic $K$-bundle $q:P\to M$ is flat if and only if it is equivalent 
to a bundle of the form $P_{\phi}$, where
$$\phi\in\Hom(\pi_1(M),K)\subset Z^1(\pi_1(M),{\mathcal O}(\widetilde{M},K)).$$
Here we identify $K$ with the subgroup of constant functions in 
${\mathcal O}(\widetilde{M},K)$.
\end{Proposition}

\begin{Proof}
The idea is to pull the bundle $P$ back to $\widetilde{M}$ and to use that a 
flat bundle on a $1$-connected manifold is trivial, see {\it loc. 
cit.} for example. The homomorphism $\phi$ is then constructed from the 
trivialization. 

On the other hand, the trivial principal $K$-bundle on $\widetilde{M}$ has the
pullback of the Maurer-Cartan form on $K$ as an invariant connection form. The
associated connection form on $M$ is the expected flat connection.
\end{Proof} 

For the following proposition, denote by $H^1_{\rm dR}(M,{\mathfrak k})$ the
set of orbits of $Z^1_{\rm dR}(M,{\mathfrak k})$ under the gauge action 
(\ref{gauge_action}) of ${\mathcal O}(M,K)$ which has been introduced in the 
previous section. The proposition appears in the smooth context in \cite{Nee04}
and in the holomorphic context as Lemma 6 in \cite{Lau}.   

\begin{Proposition} \label{bundle_exact_sequence}
Let $P$ be the map which assigns to $\alpha\in Z^1_{\rm dR}(M,{\mathfrak k})$ 
its period homomorphism, and let ${\rm Bun}_h^{\flat}(M,K)$ denote the set of 
equivalence classes of holomorphic flat principal $K$-bundles. 
Then there is an exact sequence of pointed sets
$$1\to K\to{\mathcal O}(M,K)\stackrel{\delta}{\to}Z^1_{\rm dR}(M,{\mathfrak k})
\stackrel{P}{\to}\Hom(\pi_1(M),K)\to{\rm Bun}_h^{\flat}(M,K)\to *.$$
More precisely, exactness means here that the fibers of $P$ are exactly the 
${\mathcal O}(M,K)$-orbits, and exactness in $\Hom(\pi_1(M),K)$ means that
the image of $P$ consists of those homomorphisms for which the corresponding
flat bundle is trivial. 
   
The sequence induces an exact sequence in cohomology
$$0\to H^1_{\rm dR}(M,{\mathfrak k})\stackrel{P}{\to}H^1(\pi_1(M),K)\to
{\rm Bun}_h^{\flat}(M,K)\to *.$$
\end{Proposition}

\begin{Proof}
Proposition \ref{preceeding_proposition} gives rise to the map 
$\Hom(\pi_1(M),K)\to{\rm Bun}_h^{\flat}(M,K)$. Exactness in 
$Z^1_{\rm dR}(M,{\mathfrak k})$ can be seen as follows:
Let $S:Z^1_{\rm dR}(M,{\mathfrak k})\to{\mathcal O}(\widetilde{M},K)$ the map
defined by $\delta(S(\alpha))=q_M^*\alpha$ and $S(\alpha)(\tilde{m}_0)=1$ for
some base point $\tilde{m}_0\in\widetilde{M}$ with $q_M(\tilde{m}_0)=m_0$. 
If ${\rm per}^{{m}_0}_{\alpha}=
{\rm per}^{{m}_0}_{\beta}$, then the function $S(\alpha)^{-1}S(\beta)$ is
constant on the orbits of $\pi_1(M)$, hence of the form $q_M^*f$ for some 
$f\in{\mathcal O}(M,K)$ with $f(m_0)=1$. By injectivity of $S$, one gets
$\beta=\alpha *f$, i.e. $\alpha$ and $\beta$ are in the same 
${\mathcal O}(M,K)$-orbit.

Exactness in $\Hom(\pi_1(M),K)$ follows directly from Corollary
\ref{preceeding_corollary}.
\end{Proof}

Let $\Sigma$ be a compact connected Riemann surface of genus $g$, and let
$p_1,\ldots,p_k$ be $k$ distinct points on $\Sigma$ for $k\geq 1$. We will
denote by $\Sigma_k$ the open Riemann surface 
$$\Sigma_k\,:=\,\Sigma\setminus\{p_1,\ldots,p_k\}.$$

Let us recall that $\Sigma_k$ is homotopy equivalent to a bouquet of $2g+k-1$ 
circles. The number $2g+k-1$ will be denoted $\ell$ in the following.
The first homotopy group $\pi_1(\Sigma_k)$ is thus the free group on $\ell$ 
generators. In the same way, the singular 
cohomology group $H^1(\Sigma_k,\Z)$ is the free abelian group on $\ell$ 
generators (some explanation of this will be given in the next section).   

\begin{Corollary}  \label{gauge_orbits}
For $M=\Sigma_k$ and connected $K$, the monodromy map $P$ induces a bijection
$$Z^1_{\rm dR}(\Sigma_k,{\mathfrak k})\,/\,{\mathcal O}(\Sigma_k,K)\,=:\,
H^1_{\rm dR}(\Sigma_k,{\mathfrak k})\,=\,K^{\ell}\,/\,K,$$
where $\ell$ is the first Betti number of $\Sigma_k$ and the action of 
${\mathcal O}(\Sigma_k,K)$ on $Z^1_{\rm dR}(\Sigma_k,{\mathfrak k})$ is the 
gauge action (\ref{gauge_action}) and the action of $K$ on 
$K^{\ell}$ is the conjugation in each argument.
\end{Corollary} 

\begin{Proof}
Indeed, for $M=\Sigma_k$, the above exact sequence reads
$$0\to H^1_{\rm dR}(\Sigma_k,{\mathfrak k})\stackrel{P}{\to}H^1(\pi_1(\Sigma_k)
,K)\to{\rm Bun}_h^{\flat}(\Sigma_k,K)\to 0.$$
Now on the one hand, $\pi_1(\Sigma_k)$ is the free group $F(\ell)$ on $\ell$ 
generators, where $\ell$ is the first Betti number of $\Sigma_k$,
and we have therefore
$$H^1(\pi_1(\Sigma_k),K)\,=\,\Hom(F(\ell),K)\,/\,K\,=\,K^{\ell}\,/\,K.$$
where the action of $K$ on $K^{\ell}$ is by conjugation in each argument.

On the other hand, exactness in $\Hom(\pi_1(M),K)$ means that
the image of $P$ consists of those homomorphisms for which the corresponding
flat bundle is trivial. But in our case, all holomorphic bundles are trivial,
as it is shown in the following proposition. This shows our claim.
\end{Proof}  

\begin{Proposition}  \label{all_trivial}
Let $K$ be a connected complex Banach Lie group.
Then any holomorphic principal $K$-bundle on $\Sigma_k$ is holomorphically 
trivial.
\end{Proposition}

\begin{Proof}
Let us first show that any $K$-principal bundle $P$ on 
$\Sigma_k$ is topologically trivial. The question is invariant under homotopy 
equivalence, so we may assume that the
base of the bundle is a bouquet of circles $\bigvee_{j=1}^{\ell}\bS^1$. 
Now we apply homotopy theory:
$${\rm Bun}(\bS^1,K)\,=\,[\bS^1,BK]\,=\,\pi_1(BK)\,/\,{\rm conjugation}\,=\,
\pi_0(K)\,/\,{\rm conjugation},$$
the last space clearly being a one-point space, because $K$ is connected. 
The bundle is thus trivial over each circle in the bouquet; the trivializations
may be glued together in the distinguished point to trivialize the bundle 
on the whole $\bigvee_{j=1}^{\ell}\bS^1$. 

On the other hand, by \cite{Rae} (see the details in \cite{Ram} and our 
Theorem \ref{Oka_principle}), 
any holomorphic
$K$-principal bundle on a Stein manifold is holomorphically trivial if and only 
if it is topologically trivial.
\end{Proof} 

\begin{Remark}
The first step of the above proof is actually the first step in the theory
of describing the obstructions to lifting a map from a (relative) CW complex 
$(X,A)$ to some target space by lifting it step by step from the $k$ skeleton 
of $X$ to the $k+1$ skeleton. 
It has been set up by Eilenberg in 1940, and a good
reference for this is \cite{Whi} Ch. V.5 p. 229.
\end{Remark} 

\begin{Remark}
In \cite{Oni} (see also \cite{Lau}, Lemma 7) it is stated that for a general 
complex
manifold $M$ and a general complex Lie group $K$, the period map $P$ in 
the above sequence is surjective if and only if one of the following 
conditions is holds:
\begin{itemize}
\item $\pi_1(M)$ is a free group and $K$ is connected,
\item $\pi_1(M)$ is a free abelian group and $K$ is a compact connected
group whose cohomology is torsion free. 
\end{itemize}
In particular, under these conditions every flat holomorphic principal 
$K$-bundle is holomorphically trivial. 
\end{Remark}

\section{The Lie group structure on the central extension}
  
In this section, we recall topological and geometrical preliminaries 
about the holomorphic current group ${\mathcal O}(\Sigma_k,K)$.
 
\subsection{The Lie group ${\mathcal O}(\Sigma_k,K)$} \label{Liegrp_structure}
 
Let $\Sigma$ be a compact connected Riemann surface of genus $g$, let
$p_1,\ldots,p_k$ be $k$ distinct points on $\Sigma$ for $k\geq 1$ and 
let $\Sigma_k=\Sigma\setminus\{p_1,\ldots,p_k\}$.

The open Riemann surface $\Sigma_k$ is homotopy equivalent to a bouquet 
of $2g+k-1$ circles. This follows, for example, from \cite{Han}. 
Intuitively speaking, 
the usual model for $\Sigma$ which proceeds by gluing the edges of a $4g$-gon
leads to a bouquet of circles, because taking out (at least) one point destroys
the $2$-cell. The number $2g+k-1$ is denoted $\ell$. 

The first homotopy group $\pi_1(\Sigma_k)$ is thus the free group on $\ell$ 
generators. Let us denote these generators (which we fix once and for all) 
$\alpha_1,\ldots,\alpha_{\ell}\in{\mathcal C}(\bS^1,\Sigma_k)$. 
In the same way (for example, using either the knowledge of $\pi_1$ or a 
Mayer-Vietoris sequence), the singular 
cohomology group $H^1(\Sigma_k,\Z)$ is the free abelian group on $\ell$ 
generators.  

In view of Huber's Theorem \cite{Hub} and the local contractibility of 
$\Sigma_k$, the group $H^1(\Sigma_k,\Z)$ is isomorphic to 
$\check{H}^1(\Sigma_k,\Z)\cong[\Sigma_k,\bS^1]$, the set of homotopy classes 
of continuous maps from $\Sigma_k$ to $\bS^1$. In particular, there exist 
continuous maps $f_1,\ldots,f_{\ell}:\Sigma_k\to\bS^1$ such that 
$[f_j\circ\alpha_i]=\delta_{ij}\in\pi_1(\bS^1)\cong\Z$. As every homotopy class
$[\Sigma_k,\bS^1]$ contains a smooth function (which may be seen by using a homotopy equivalence of $\Sigma_k$ with a compact surface with boundary),
we may assume in the following that the $f_j$ are smooth. As the 2-pointed 
Riemann sphere $\bS^2\setminus\{0,\infty\}\,=\,\C^{\times}$ is homotopy 
equivalent to the circle $\bS^1$, 
$$[\Sigma_k,\bS^1]\,=\,[\Sigma_k,\C^{\times}],$$
and each homotopy class contains also a holomorphic function (see \cite{Rae}
and Theorem \ref{Oka_principle}).
In this sense, we may also assume the functions $f_j$ to be holomorphic.  

Let $Y$ be a sequentially complete locally convex space.    
Using the logarithmic derivatives $\delta(f_j)$, seen here as usual closed 
holomorphic $1$-forms on $\Sigma_k$,
one obtains an isomorphism
$$\Phi:H_{\rm dR}^1(\Sigma_k,Y)\to Y^{\ell},\,\,\,\,\,[\beta]\mapsto
\left(\int_{\alpha_j}
\beta\right)_{j=1,\ldots,\ell},$$ 
with explicit inverse $\Phi^{-1}(y_1,\ldots,y_{\ell})=\left[\sum_{j=1}^{\ell}
\delta(f_j)y_j\right]$.

In the same way as the compact differentiable manifold $\Sigma$ may also be
regarded as a complete projective algebraic curve over $\C$, $\Sigma_k$ may 
be regarded as an affine algebraic curve over $\C$. Therefore $\Sigma_k$ 
carries a natural structure of an algebraic variety and of a complex manifold,
and it makes sense to associate to $\Sigma_k$ its structure sheaves of 
regular functions ${\rm Reg}_{\Sigma_k}$ and of holomorphic functions 
${\mathcal O}_{\Sigma_k}$, and their vector spaces of global sections
${\rm Reg}_{\Sigma_k}(\Sigma_k)$ and ${\mathcal O}_{\Sigma_k}(\Sigma_k)$,
or ${\rm Reg}(\Sigma_k)$ and ${\mathcal O}(\Sigma_k)$ in short hand 
notation.
Both of them are in general infinite dimensional vector spaces over $\C$.

Let $K$ be a connected complex Banach Lie group, and ${\mathfrak k}$ be
its Lie algebra. Most of the time, we will suppose that $K$ is 1-connected,
finite-dimensional
and that ${\mathfrak k}$ is {\it simple}, i.e. does not possess a 
non-trivial proper ideal.
  
\begin{Definition}
The holomorphic current group ${\mathcal O}(\Sigma_k,K)$ is the topological 
group of holomorphic maps on $\Sigma_k$ with values in $K$. The group structure
is given by the pointwise multiplication. The topology is 
given by the compact open topology on the mapping space, or equivalently by the
topology of uniform convergence on compact subsets. We will consider 
the holomorphic current algebra ${\mathcal O}(\Sigma_k,{\mathfrak k})$ as 
the Lie algebra of ${\mathcal O}(\Sigma_k,K)$. The bracket on 
${\mathcal O}(\Sigma_k,{\mathfrak k})$ is given pointwise. 
${\mathcal O}(\Sigma_k,K)$ together with the compact open topology is a 
Fr\'echet-Lie algebra.
\end{Definition}

In fact, ${\mathcal O}(\Sigma_k,K)$ is a Lie group, i.e. carries a Fr\'echet
manifold structure such that multiplication and inversion maps are smooth. 
This is the content of the following theorem which is Theorem $3.12$ in 
\cite{NeeWag}. It is also shown in {\it loc. cit.}
that ${\mathcal O}(\Sigma_k,K)$ is a {\it regular Lie group}, i.e. that the 
evolution equation (see Section 1) always has solutions, i.e. the evolution
map is well defined.  

Observe also that 
${\mathcal O}(\Sigma_k,K)$ carries an exponential map 
$$\exp:{\mathcal O}(\Sigma_k,{\mathfrak k})\,\to\,{\mathcal O}(\Sigma_k,K),$$
which is just given by the concatenation with the finite dimensional 
exponential map $\exp_K:{\mathfrak k}\to K$.

Let now $K$ be a connected affine algebraic group over $\C$.  
It seems reasonable to expect that ${\rm Reg}(\Sigma_k,K)$ 
also carries a structure
of an infinite-dimensional algebraic group in case $K$ is an algebraic group. 
One approach in this direction would describe it 
as a group valued functor ${\tt rings}\to{\tt groups}$, sending a ring $R$
to the group ${\rm Reg}(\Sigma_k(R),K))$ of regular maps from the $R$-points of 
the affine algebraic curve $\Sigma_k$ to the algebraic group $K$, 
but no results in 
this direction are known to us. 

\begin{Theorem}
Let $M$ be a non-compact connected complex curve without boundary. Assume 
further that $\pi_1(M)$ is finitely generated and that $K$ is a complex 
Banach-Lie group. Then the group ${\mathcal O}_*(M,K)$ carries a 
regular complex Lie group structure for which
$$\delta:{\mathcal O}_*(M,K)\to \Omega^1(M,{\mathfrak k})$$
is biholomorphic onto a complex submanifold, and 
${\mathcal O}(M,K)\cong K\ltimes{\mathcal O}_*(M,K)$ carries
a regular complex Lie group structure compatible with evaluations and the 
compact open topology.
\end{Theorem}

Of course, we apply this theorem here with $M=\Sigma_k$. 

The manifold structure on the space of pointed maps
${\mathcal O}_*(M,K)$ is obtained using an infinite dimensional version
of the regular value theorem (using a parametrized version of the implicit 
function theorem of Gl\"ockner \cite{Glo}). 
More precisely, the image of the logarithmic derivative 
$$\delta:{\mathcal O}_*(M,K)\to \Omega^1(M,{\mathfrak k})$$
in the Fr\'echet space of holomorphic $1$-forms $\Omega^1(M,{\mathfrak k})$
with values in ${\mathfrak k}$ is characterized as the inverse image of the 
trivial monodromy homomorphism under the period map, and is thus seen to be a 
split submanifold.

\subsection{Topology of the Lie group ${\mathcal O}(\Sigma_k,K)$}

Here we compute the homotopy groups of the Lie group ${\mathcal O}(\Sigma_k,K)$.
It is well-known from Sullivan's rational homotopy theory that the rational 
cohomology algebra of the H-space
${\mathcal O}(\Sigma_k,K)$ has as its generators the duals of the generators of
the $\pi_i({\mathcal O}(\Sigma_k,K))\otimes\Q$, thus these computations 
determine in particular the rational cohomology algebra. 

All our computations are based on the Oka principle (see \cite{Rae}, with 
details of the proof from \cite{Ram}):

\begin{Theorem}  \label{Oka_principle}
Let $M$ be a Stein manifold and $K$ be a connected complex Banach Lie group.
Then
\begin{enumerate}
\item[(a)] each continuous map $f:M\to K$ is homotopic to a holomorphic map;
\item[(b)] if two holomorphic maps $f_0,f_1:M\to K$ are homotopic within 
continuous maps (i.e. in ${\mathcal C}^0(M,K)$), then they are homotopic within
holomorphic maps (i.e. in ${\mathcal O}(\Sigma_k,K)$). 
\end{enumerate}
\end{Theorem}

\begin{Corollary}
Let $M$ be a Stein manifold and $K$ be a connected complex Banach Lie group. 
Then ${\mathcal O}(\Sigma_k,K)$ has the same homotopy type as 
${\mathcal C}^0(M,K)$.
\end{Corollary}

As before, we also use for the computation of 
$\pi_i({\mathcal O}(\Sigma_k,K))$ that $\Sigma_k$ is homotopy equivalent to 
a bouquet of $\ell=2g+k-1$ circles $\bS^1$. The outcome is then:

\begin{Proposition}
Let $K$ be a connected complex Banach Lie group.
Then for $i\geq 0$,
$$\pi_i({\mathcal O}(\Sigma_k,K))\,=\,\pi_i(K)\oplus\pi_{i+1}(K)^{\ell}.$$
\end{Proposition}

\begin{Proof}
The group ${\mathcal O}(\Sigma_k,K)$ has the same homotopy type as 
${\mathcal C}^0(\Sigma_k,K)$. The topological group 
${\mathcal C}^0(\Sigma_k,K)$ splits as
a semi-direct product 
$${\mathcal C}^0(\Sigma_k,K)\,=\,{\mathcal C}^0_*(\Sigma_k,K)\rtimes K,$$
where ${\mathcal C}^0_*(\Sigma_k,K)$ denotes the space of continuous maps 
$f:\Sigma_k\to K$ such that $f(*)=1$ for the base point $*$ of $\Sigma_k$.
 
Now compute (with $S$ being the suspension and $[,]$ the based homotopy classes)
\begin{eqnarray*}
\pi_i({\mathcal C}^0_*(\Sigma_k,K))&=&[S^i(\vee_{k=1}^{\ell}\bS^1),K]  \\
&=&[\bS^{i+1}\vee\ldots\vee\bS^{i+1},K] \\
&=& \pi_{i+1}(K)^{\ell}.
\end{eqnarray*}

This ends the proof.
\end{Proof}
  
\begin{Corollary} \label{1-connected}
Let $M$ be a Stein manifold and $K$ be a connected complex Banach Lie group.
Then ${\mathcal O}(\Sigma_k,K)$ is 1-connected if and only if $K$ is 
2-connected. This is the case, for example, if $K$ is finite dimensional
and 1-connected. 
\end{Corollary}

\subsection{Central extensions of current algebras}

Here we recall the second cohomology of current algebras. The earliest 
reference is usually attributed to Bloch \cite{Blo} and Kassel-Loday 
\cite{Kas}. 
While these references deal with homology, for 
cohomology, and especially continuous cohomology, we refer to \cite{NeeWag1}.

Let $A$ be a unital commutative associative algebra over a field $\K$ of 
characteristic zero and ${\mathfrak k}$ a $\K$-Lie algebra. The tensor product 
${\mathfrak g}:=A\otimes{\mathfrak k}$ is a Lie algebra when endowed with the
{\it current bracket}
$$[a\otimes x,b\otimes y]\,=\,ab\otimes[x,y].$$
Let ${\mathfrak g}'$ denote the derived algebra of ${\mathfrak g}$, i.e. the 
commutator ideal.
The goal of this subsection is to present the second cohomology space 
$H^2({\mathfrak g})$ in terms of data associated to $A$ and ${\mathfrak k}$.

For this purpose, let $\Omega^1(A)$ be the {\it module of K\"ahler 
differentials} associated to $A$. It may be defined by its universal property,
namely for the $A$-module $\Omega^1(A)$, there is a derivation 
$d_A:A\to\Omega^1(A)$ such that for any other $A$-module $M$ and any other
derivation $D:A\to M$, there is a unique morphism of $A$-modules 
$\alpha:\Omega^1(A)\to M$ such that $D=\alpha\circ d_A$. There are several
explicit constructions of $\Omega^1(A)$, and we refer to \cite{NeeWag1} for 
the most commonly known. For a Fr\'echet algebra $A$, there is a version of
$\Omega^1(A)$ which takes the topology into account and is a Fr\'echet 
$A$-module, see {\it loc. cit.} for example.

For the Lie algebra ${\mathfrak k}$, we have the usual cycle, boundary and 
cohomology spaces $Z^n(\k)$, $B^n(\k)$ and $H^n(\k)=Z^n(\k)\,/\,B^n(\k)$, all
with trivial coefficients. Recall further from \cite{NeeWag1} that there is
a map, called the {\it Cartan map}
$$\Gamma:{{\rm Sym}}({\mathfrak k})^{\mathfrak k}\to Z^3(
{\mathfrak k}),\,\,\,\,\,\,\,\Gamma(\kappa)(x,y,z)\,:=\,\kappa([x,y],z),$$
linking the space ${\rm Sym}({\mathfrak k})^{\mathfrak k}$ of 
${\mathfrak k}$-invariant symmetric bilinear forms to the space of 
$3$-cocycles. Denote by $Z^3(\k)_\Gamma$ and $B^3(\k)_\Gamma$ the intersections
$Z^3(\k)_\Gamma:=Z^3(\k)\cap\im(\Gamma)$ and 
$B^3(\k)_\Gamma:=B^3(\k)\cap\im(\Gamma)$ with the image of $\Gamma$.

One of the main results of \cite{NeeWag1} is the following Theorem 
(Theorem $4.2$):

\begin{theorem}  The sequence 
$$ 0 \to H^2(\g/\g') \oplus \big(A \otimes H^2(\k)) \ssmapright{} H^2(\g) 
\ssmapright{} \Lin((\Omega^1(A), d_A(A)), (Z^3(\k)_\Gamma, B^3(\k)_\Gamma))
\to 0$$ 
is exact. 
\end{theorem}

In the special case where ${\mathfrak k}$ is a finite dimensional complex
simple Lie algebra, the space ${{\rm Sym}}({\mathfrak k})^{\mathfrak k}$ is
$1$-dimensional, ${\mathfrak k}'={\mathfrak k}$ and $H^2({\mathfrak k})=0$,
therefore we get  
$$H^2(\g)\,=\,(\Omega^1(A)\,/\,d_A(A))^*,$$
cf \cite{Blo} and \cite{Kas}.

The generating cocycle for the universal central extension is described as 
follows. Let $\kappa$ be the Cartan-Killing form on ${\mathfrak k}$. Then 
the $2$-cocycle $\omega_{\kappa}$ is defined by
$$
\omega_{\kappa}(x,y)\,=\,[\kappa(x,d_Ay)],
$$
where $[.]$ denotes the equivalence class in $\Omega^1(A)\,/\,d_A(A)$.

We will apply this to the case where the associative algebra $A$ is the 
Fr\'echet algebra
of holomorphic functions ${\mathcal O}(\Sigma_k)$ on the punctured Riemann 
surface
$\Sigma_k$. In this case, the Fr\'echet space of K\"ahler differentials can be 
identified with the Fr\'echet space of differential $1$-forms. This is the 
content of the main theorem in \cite{NeeWag3} (Theorem $2.1$):  

\begin{Theorem}
Let $X$ be a Stein manifold. Then the map
$$\gamma_X:\Omega^1({\mathcal O}(X))\to\Omega^1(X),$$
induced by the universal property, is an isomorphism of Fr\'echet
${\mathcal O}(X)$-modules.
\end{Theorem} 

Let ${\mathfrak k}$ be a finite dimensional simple complex Lie algebra.
The main object of the present article is the central extension
$\widehat{{\mathcal O}(\Sigma_k,{\mathfrak k})}$ of 
${\mathcal O}(\Sigma_k,{\mathfrak k})$ 
by means of the cocycle $\omega_{\kappa}$. In the next subsection, we 
will describe how this central extension is integrated into a Fr\'echet-Lie 
group $\widehat{{\mathcal O}(\Sigma_k,K)}$, central extension of the Lie group 
${\mathcal O}(\Sigma_k,K)$.

\subsection{The Lie group structure on the central extension}
\label{section_central_extension}

In this subsection, we will adapt the Reduction Theorem of \cite{NeeMai} to 
our holomorphic framework. This will imply that the Lie group structure on
${\mathcal O}(\Sigma_k,K)$ which was described in subsection 
\ref{Liegrp_structure} extends to the central extension. 

Let $M$ be a complex manifold, $K$ be a complex Lie group and $Y$ be a
sequentially complete locally convex complex space. 
We will work with the topological group $G:={\mathcal O}(M,K)$ of holomorphic 
maps from $M$ to $K$ and its Lie algebra 
${\mathfrak g}:={\mathcal O}(M,{\mathfrak k})$.
Let $\kappa:{\mathfrak k}\times{\mathfrak k}\to Y$ be a continuous invariant 
symmetric bilinear $Y$-valued form on ${\mathfrak k}$.    
The cocycles which we consider on ${\mathfrak g}$ are of the form
\begin{equation}     \label{cocycle}
\omega_{M}(\xi,\eta)\,:=\,
[\kappa(\xi,d\eta)]\,\in\,\Omega^1(M,Y)\,/\,d\,{\mathcal O}(M,Y),
\end{equation}
where $\xi,\eta\in{\mathfrak g}$, and $[.]$
denotes the equivalence class in the quotient space 
$${\mathfrak z}_M(Y)\,:=\Omega^1(M,Y)\,/\,d\,{\mathcal O}(M,Y).$$ 
A typical example for $\kappa$ is the Cartan-Killing form on a simple complex 
Lie algebra ${\mathfrak k}$, where $Y=\C$. 
In the discussion that follows, cocycles of the same kind which are defined on
mapping spaces (of smooth maps) with domain the circle $\bS^1$ will play a 
role. We will denote them accordingly $\omega_{\bS^1}$. 

For the following lemma, observe that the usual loop group
${\mathcal C}^{\infty}(\bS^1,K)$ is a {\it complex} Fr\'echet-Lie group in our 
case, because $K$ is a complex Lie group. Recall that a map $f:X\to Y$ 
between complex Fr\'echet manifolds is called 
{\it holomorphic} in case $f$ is smooth and its differentials are complex 
linear.

\begin{Lemma} \label{homomorphism}
Let $\alpha:\bS^1\to M$ be smooth.  
 
The group homomorphism $\alpha_K:G\to{\mathcal C}^{\infty}(\bS^1,K)$ defined by
$f\mapsto f\circ\alpha$ is a holomorphic homomorphism of complex Lie groups, 
where $G$ carries the Fr\'echet-Lie group structure described in Section 
\ref{Liegrp_structure} and 
${\mathcal C}^{\infty}(\bS^1,K)$ carries the usual (loop group) Fr\'echet-Lie 
group structure.
\end{Lemma} 

\begin{Proof}
Observe that $\alpha_K$ is the restriction map, i.e. the restriction 
$f|_{\alpha(\bS^1)}$ of 
holomorphic functions $f\in G={\mathcal O}(M,K)$ to the ``circle'' 
$\alpha(\bS^1)$. 
As such, it is first of all continuous, because if $f_i\to f$ 
in the compact-open topology in $G$, then $(f_i)^{[n]}\to f^{[n]}$ in the 
smooth compact-open topology. But then this is also true for all restrictions 
$(f_i)^{[n]}|_{\alpha(\bS^1)}\to f^{[n]}|_{\alpha(\bS^1)}$, and therefore 
$\alpha_K$ is continuous.

Let us prove that $\alpha_K$ is smooth. This actually follows from the same
restriction argument as above, because the tangent map
$$T\alpha_K:T{\mathcal O}(M,K)\to T{\mathcal C}^{\infty}(\bS^1,K)$$
can be identified with a two component restriction map under the isomorphism
$$T{\mathcal O}(M,K)\,\cong\,{\mathcal O}(M,K)\times{\mathcal O}(M,
{\mathfrak k})$$
and in the same way
$$T{\mathcal C}^{\infty}(\bS^1,K)\,\cong\,{\mathcal C}^{\infty}(\bS^1,K)\times
{\mathcal C}^{\infty}(\bS^1,{\mathfrak k}).$$
Iteration of this argument implies that $\alpha_K$ is smooth.

Now a map $f:X\to Y$ between complex Fr\'echet manifolds is holomorphic
in case $f$ is smooth and its differentials are complex 
linear. As the restriction map is complex linear on the tangent spaces, our
map $\alpha_K$ is therefore holomorphic. 
\end{Proof} 

\begin{Remark}
For $M=\Sigma_k$ and a smooth map $\alpha:\bS^1\to M$ whose image contains an 
interval which is not reduced to a point, the Lie group homomorphism 
$\alpha_K:G\to{\mathcal C}^{\infty}(\bS^1,K)$ is injective. Indeed, in case the
restriction of two holomorphic functions to the set $\alpha(\bS^1)$ coincide, 
the two functions coincide, because the set $\alpha(\bS^1)$ contains an 
accumulation point. 
\end{Remark} 

By definition, $\omega_{M}$ is a ($Y$-valued) $2$-form on the tangent space of 
$G={\mathcal O}(M,K)$ at $1\in G$, and may thus be extended to a unique 
left-invariant ($Y$-valued) $2$-form $\Omega_{M}$ on $G$ such that 
$\Omega_{M}(1)=\omega_{M}$.

Choosing in a homotopy class a smooth representative permits to define a 
{\it period homomorphism}
$${\rm per}_{\omega_{M}}:\pi_2(G)\to{\mathfrak z}_M(Y),\,\,\,\,
{\rm per}_{\omega_{M}}([\sigma])\,:=\,\int_{\sigma}\Omega_{M}.$$
Observe that the homotopy group $\pi_2(G)$ of the infinite dimensional Lie 
group $G$ is not necessarily zero. In order to extend a given Lie group 
structure on ${\mathcal O}(M,K)$ to the central extension given by the cocycle
$\omega_{M}$, it is necessary and sufficient that the {\it period group}
$\Pi_{M,\kappa}:=\im({\rm per}_{\omega_{M}})$, i.e. the image of the 
period homomorphism ${\rm per}_{\omega_{M}}$, is a discrete subspace in 
${\mathfrak z}:={\mathfrak z}_M(Y)$, see \cite{Neeb}. 
If this condition is fulfilled, 
the central extension is then 
an extension of ${\mathcal O}(M,K)$ by the abelian group 
${\mathfrak z}\,/\,\Pi_{M,\kappa}$. 

We now adapt the discussion of Section I in \cite{NeeMai} to our context.
For $\alpha\in{\mathcal C}^{\infty}(\bS^1,M)$, denote by 
$\alpha_{\mathfrak z}:{\mathfrak z}\to\C$ the integration over 
$\alpha(\bS^1)$.

\begin{Lemma}
For each $\alpha\in{\mathcal C}^{\infty}(\bS^1,M)$, we have 
$$\alpha_{\mathfrak z}\circ{\rm per}_{\omega_{M}}\,=\,
{\rm per}_{\omega_{\bS^1}}\circ\pi_2(\alpha_K),$$
where $\pi_2(\alpha_K):\pi_2(G)\to\pi_2({\mathcal C}^{\infty}(\bS^1,K))$ is the
group homomorphism induced by the Lie group homomorphism 
$\alpha_K:G\to{\mathcal C}^{\infty}(\bS^1,K)$. 
\end{Lemma}

\begin{Proof}
The proof of Lemma I.7 in \cite{NeeMai} is easily adapted: 
For $\alpha\in{\mathcal C}^{\infty}(\bS^1,M)$, the map
$\alpha_{\mathfrak z}:{\mathfrak z}\to\C$ denotes the integration over 
$\alpha(\bS^1)$. Then $\alpha_{\mathfrak z}\circ\Omega_{M}$ is a left-invariant
$2$-form on $G$ whose value at $1\in G$ is 
$\alpha_{\mathfrak z}\circ\omega_{M}$. Further 
$\alpha_K^*\Omega_{\bS^1}$ is a left-invariant $2$-form on $G$ whose value at 
$1\in G$ is given by
\begin{eqnarray*}
(\xi,\eta)\mapsto\omega_{\bS^1}(\xi\circ\alpha,\eta\circ\alpha)&=&
[\kappa(\xi\circ\alpha,d(\eta\circ\alpha))]  \\
&=& [\kappa(\alpha^*\xi,\alpha^*d\eta)]\,=\,\int_{\bS^1}\kappa
(\alpha^*\xi,\alpha^*d\eta)\\
&=&\int_{\alpha(\bS^1)}\kappa(\xi,d\eta)\,=\,\alpha_{\mathfrak z}(\omega_{M}
(\xi,\eta)).
\end{eqnarray*}
This implies $\alpha_{\mathfrak z}\circ\Omega_{M}=\alpha^*_K
\Omega_{\bS^1}$, which shows the lemma.
\end{Proof} 

Let us now restrict to $\Sigma_k$.
Following the reasoning of Section I in \cite{NeeMai}, Corollary I.9 remains 
true for $\Sigma_k$, because holomorphic $1$-forms are closed. 
Lemma I.10 does not involve the holomorphic context, thus we get:

\begin{theor}[Reduction Theorem]
Identifying $H^1_{\rm dR}(\Sigma_k,Y)$ with $Y^{\ell}$ via 
the isomorphism $\Phi$, where $\ell$ is the first Betti number of $\Sigma_k$, 
we have
$$\Pi_{\Sigma_k,\kappa}\,\cong\,\Pi_{\bS^1,\kappa}^{\ell}\subset Y^{\ell}\,\cong\,
H^1_{\rm dR}(\Sigma_k,Y)\,=\,{\mathfrak z}_{\Sigma_k}(Y).$$
In particular, $\Pi_{\Sigma_k,\kappa}$ is discrete if and only if 
$\Pi_{\bS^1,\kappa}$ is discrete. 
\end{theor}

\begin{Proof}
To each $f\in{\mathcal C}^{\infty}(\Sigma_k,\bS^1)$, one associates the map
$f_K:{\mathcal C}^{\infty}(\bS^1,K)\to G={\mathcal O}(\Sigma_k,K)$, $\eta\mapsto 
\eta\circ f$, which in turn induces a map $\pi_2(f_K):\pi_2(
{\mathcal C}^{\infty}(\bS^1,K))\to \pi_2(G)$. For $\alpha\in
{\mathcal C}^{\infty}(\bS^1,\Sigma_k)$, we obtain with Lemma I.10 of \cite{NeeMai}
\begin{eqnarray*}
\alpha_{\mathfrak z}\circ {\rm per}_{\omega_{\Sigma_k}}\circ\pi_2(f_K)&=&
{\rm per}_{\omega_{\bS^1}}\circ\pi_2(\alpha_K)\circ\pi_2(f_K)=
{\rm per}_{\omega_{\bS^1}}\circ\pi_2(\alpha_K\circ f_K)\\
&=& {\rm per}_{\omega_{\bS^1}}\circ\pi_2((f\circ\alpha)_K)={\rm deg}
(f\circ\alpha)\,{\rm per}_{\omega_{\bS^1}}
\end{eqnarray*}
 
For $f=f_j$ and $\alpha=\alpha_i$, it follows in particular that 
$\alpha_{i,{\mathfrak z}}\circ{\rm per}_{\omega_{\Sigma_k}}\circ\pi_2(f_{j,K})=
\delta_{ij}\,{\rm per}_{\omega_{\bS^1}}$. Here $\alpha_{i,{\mathfrak z}}$ may be 
seen as the projection onto the $i$-th factor in 
$$Y^{\ell}\,\cong\,H^1_{\rm dR}(\Sigma_k,Y)=\bigoplus_{j=1}^{\ell}
[\delta(f_j)]\,Y.$$ 
Hence
$${\rm per}_{\omega_{\Sigma_k}}(\im\,\pi_2(f_{j,K}))\,=\,[\delta(f_j)]\,\Pi_{\bS^1}$$
and further $\Pi_{\Sigma_k}\supset\sum_{j=1}^{\ell}[\delta_{ij}]\Pi_{\bS^1}=
\Pi_{\bS^1}^{\ell}$.

For the reverse inclusion, we observe that $\alpha_{j,{\mathfrak z}}\circ 
{\rm per}_{\omega_{\Sigma_k}}={\rm per}_{\omega_{\bS^1}}\circ\pi_2(\alpha_K)$ implies
that for each $j$, we have $\alpha_{j,{\mathfrak z}}\circ {\rm per}_{\omega_{\Sigma_k}}
\subset \Pi_{\bS^1}$, and therefore $\Pi_{{\Sigma_k}}\subset\Pi_{\bS^1}^l$. 
\end{Proof} 

\begin{Corollary}
Suppose now that the complex Lie group $K$ is 1-connected and finite 
dimensional. Then $\Pi_{\Sigma_k,\kappa}$ is discrete, and the Fr\'echet-Lie 
group structure on
${\mathcal O}(\Sigma_k,K)$ lifts to a Lie group structure on the central 
extension.
\end{Corollary} 

\begin{Proof}
Indeed, it is shown in Theorem II.9 in \cite{NeeMai} that $\Pi_{\bS^1,\kappa}$ 
is discrete. For the second part of the corollary, observe that under the 
above hypotheses, ${\mathcal O}(\Sigma_k,K)$ is 1-connected by Corollary
\ref{1-connected}. It is explained in Section 6 of \cite{Neeb} how to construct
out of a Lie algebra 2-cocycle on the Lie algebra of a 1-connected Lie group 
a locally smooth group 2-cocycle with values in
$Z:={\mathfrak z}_{\Sigma_k}(Y)\,/\,\Pi_{\Sigma_k,\kappa}$. 
The corresponding central extension is
then a Lie group, as ${\mathcal O}(\Sigma_k,K)$ and the quotient 
$Z$ are Lie groups.    
\end{Proof}

As stated before, we will denote by $\widehat{{\mathcal O}(\Sigma_k,K)}$ the 
central group extension thus obtained (using the cocycle (\ref{cocycle})). 
The Lie group 
$\widehat{{\mathcal O}(\Sigma_k,K)}$ is a regular Fr\'echet-Lie group. 

\begin{Remark}
\begin{enumerate}
\item[(a)] We believe that the same reasoning shows that Lie group structures 
on ${\mathcal O}(M,K)$ for higher dimensional $M$ extend to central extensions 
by cocycles of the above form.
It would be interesting to extract from the discussion in Section $4$ in 
\cite{NeeWag} an interesting example of a Lie group ${\mathcal O}(M,K)$ with
$M$ of dimension greater than $1$ (for example with $K$ solvable), and from 
\cite{NeeWag1} a class of cocycles
on its Lie algebra ${\mathcal O}(M,{\mathfrak k})$, where a higher dimensional
complex version of the Reduction Theorem applies.
\item[(b)] It is probably easy to adapt the results of Section IV of 
\cite{NeeMai} to show that the central extension 
$\widehat{{\mathcal O}(\Sigma_k,K)}$ of
${\mathcal O}(\Sigma_k,K)$ via the universal cocycle (\ref{cocycle})
is in fact the universal central extension of ${\mathcal O}(\Sigma_k,K)$.
\end{enumerate}  
\end{Remark}

\section{Automorphism group of $\widehat{{\mathcal O}(\Sigma_k,K)}$}

In this section, we study the automorphisms of 
$\widehat{{\mathcal O}(\Sigma_k,{\mathfrak k})}$ and 
$\widehat{{\mathcal O}(\Sigma_k,K)}$. The discussion uses key results from
 \cite{NeeMai} Section V.

Let us first of all recall the automorphism group of the open Riemann surface 
$\Sigma_k=\Sigma\setminus\{p_1,\ldots,p_k\}$, obtained from a compact connected
Riemann surface $\Sigma$ by extracting $k\geq 1$ distinct points 
$p_1,\ldots,p_k$.

\begin{prop}
$\Aut(\Sigma_k)$ is a subgroup of $\Aut(\Sigma)$, namely the subgroup of 
automorphisms fixing all the points $p_1,\ldots,p_k$.
\end{prop}

\begin{Proof}
Let $\phi:\Sigma_k\to\Sigma_k$ be an automorphism, i.e. a biholomorphic 
bijective map. Let us show that $\phi$ extends to an automorphism 
$\hat{\phi}:\Sigma\to\Sigma$. 
It suffices to consider one extracted point $p\in\Sigma$, i.e. 
$k=1$. Denote by $D^*$ a punctured disc, centered in $p=0$, around $p$.
Also, denote by $D_{\epsilon}^*$ a punctured disc of radius $\epsilon$. 

The compact Riemann surface $\Sigma$, viewed as a projective curve, may be 
embedded into some projective space $\P^n$ using 
$i:\Sigma\hookrightarrow\P^n$. Take a minimal embedding, i.e. 
$n$ minimal. Denote by $[z_0:\ldots:z_{n}]$ the projective coordinates on 
$\P^n$, and observe that $u_{ij}=\frac{z_i}{z_j}:\P^n\to\P^1$ gives a regular 
map.    

Composing $\phi$ with $u_{ij}$, $\phi$ gives rise to a map 
$\phi_{ij}:D^*\to\P^1$ with a possible singularity in the point $p=0$. We have 
to show that $p=0$ is non singular for $\phi_{ij}$. 

Let us first note that $\phi_{ij}$ cannot have an essential singularity in 
$p=0$. Indeed, if a holomorphic function $f:D^*\to\C$ has an essential 
singularity in $0$, then for all $\epsilon>0$, $f(D_{\epsilon}^*)\subset\C$ is 
dense. (If this was not the case, let $c\in\C$ be not hit by $f$. Then 
$$g(z):=\frac{1}{f(z)-c}$$
is bounded on $D_{\epsilon}^*$, and therefore $f(z)=\frac{1}{g(z)}+c$ 
meromorphic, which is a contradiction.) This density is incompatible with the
bijective character of $\phi$.  

Now consider $\phi$ as the composition
$$D^*\stackrel{\phi}{\hookrightarrow}\Sigma\hookrightarrow\P^n.$$
Composed with $u_{ij}$ we obtain a finite covering 
$$v_{ij}:=u_{ij}\circ\phi:D^*\to\P^1.$$
Let ${\rm ord}_0(v_{10}),\ldots,{\rm ord}_0(v_{n0})$ be the orders (of poles 
or zeros) of the function $v_{i0}$ in $p=0\in D^*$. By what we have already 
shown, there exists $i$ such that  ${\rm ord}_0(v_{i0})$ is maximal. Define
then $w_{ij}:=\frac{v_{i0}}{v_{j0}}$. We have by construction 
${\rm ord}_0(w_{ij})\geq 0$ for all $i,j$. Therefore the $w_{ij}$ give affine
coordinates (given by holomorphic functions on $D^*$ without pole in $0$)
such that $\phi$ factors as 
$$D^*\to U\hookrightarrow\P^n$$
for some affine open set $U$. This shows that $\phi_{ij}$ does not have a pole 
in $p=0$ and can thus be extended to a holomorphic function on the entire
disc $D$. 

It is clear that this extension $\phi$ has to fix the points $p_i$
one by one.   
\end{Proof}

\begin{Corollary}
For genus $g>1$ and $g=1$, the automorphism group $\Aut(\Sigma_k)$ is finite. 
For $g=0$, $\Aut(\Sigma_k)$ is finite as soon as $k>2$.  
\end{Corollary}

\begin{Proof}
It is the content of Theorem 2.5, p. 88, in \cite{Kob} that $\Aut(\Sigma)$
is a finite group for genus $g\geq 2$. As the automorphisms in $\Aut(\Sigma_k)$ 
must fix the distinguished points $p_1,\ldots,p_k\in\Sigma$ and $k\geq 1$, 
the automorphism group of $\Sigma_k$ in genus $g=1$ must be finite, because 
$\Aut(\Sigma)\cong\C$ for $g=1$. For genus $g=0$, 
$\Aut(\Sigma)\cong{\rm PSl}(2,\C)$ is 3 dimensional and an element of 
${\rm PSl}(2,\C)$ 
fixing 3 points must be the identity. 
\end{Proof}

\begin{Proposition}
Let ${\mathfrak k}$ be a finite dimensional, semi-simple complex Lie algebra and
$M$ be a complex Stein space.
 
Then
$$\Aut({\mathcal O}(M,{\mathfrak k}))\cong \Aut(M)
\ltimes{\mathcal O}(M,\Aut({\mathfrak k})).$$
\end{Proposition}

\begin{Proof}
Let $\alpha:{\cal O}(M,{\mathfrak k})\to{\cal O}(M,{\mathfrak k})$ 
be an automorphism. 
Now consider the composition
$$a_x:{\mathfrak k}\hookrightarrow {\cal O}(M,{\mathfrak k})
\stackrel{\alpha}{\to}{\cal O}(M,{\mathfrak k})
\stackrel{{\rm ev}_x}{\to}{\mathfrak k},$$
where ${\rm ev}_x:{\cal O}(M,{\mathfrak k})\to {\mathfrak k}$ denotes the 
evaluation in $x\in M$.
The map $a_x$ is not zero, because $a_x=0$ implies 
${\rm ev}_x\circ\alpha=:\alpha_x=0$
by using ${\cal O}(M,{\mathfrak k})=[{\mathfrak k},
{\cal O}(M,{\mathfrak k})]$ (see Lemma \ref{lemma1} below).
Replacing ${\mathfrak k}$ by a suitable subalgebra, this shows that 
$a_x$ is an automorphism. 

Clearly, $x\mapsto a_x$ is an element $a$ of ${\cal O}(M,{\rm Aut}
({\mathfrak k}))$, and 
replacing $\alpha$ by $a^{-1}\circ\alpha$, one can assume 
$a_x={\rm id}_{\mathfrak k}$ for all $x\in M$. 

Let us show now that $\alpha_x={\rm ev}_y$ for some $y\in M$. 
$\alpha_x$ is a homomorphism of Lie algebras.
Its kernel must be of the form $I\otimes\fk$ where $I$ is a maximal ideal of
$A$, by Lemma \ref{lemma2} below. Therefore, 
by the theory of Stein algebras (the ``Versch\"arfung von Satz 2'', p. 181, 
in \cite{GraRem} states that   
a closed ideal is necessarily of the form $\fm_y$, but on the other hand, 
an algebra homomorphism $\pi:{\cal O}(M)\to\C$ is automatically continuous,
cf \cite{GraRem} p. 187) $I=\fm_y$ for some $y\in M$,
where $\fm_y=\{f\in{\cal O}(M)\,|\,f(y)=0\}$. It is then clear that 
$\alpha_x={\rm ev}_y$, and that ${\alpha_x}={\rm ev}_y$

We define now $\psi:M\to M$ by $x\mapsto y$ where $y$ is such that 
${\alpha_x}={\rm ev}_y$. With this definition, $\alpha$ is written as
$$\alpha(f)(x)\,=\,f(\psi(x))\,=\,\psi^*(f)(x).$$
Now $\psi$ must be holomorphic, as it sends holomorphic maps to holomorphic 
maps (this follows easily using projections as holomorphic functions), and 
finally, $\psi$ must be an automorphism of $M$. 
\end{Proof}

\begin{Lemma}   \label{lemma1}
${\cal O}(M,\fk)=[\fk,{\cal O}(M,\fk)]$
\end{Lemma}

\begin{Proof}
This follows directly from $\fk=[\fk,\fk]$.
\end{Proof}

The following lemma is very close to Lemma 6.1 in \cite{Nee}: 
 
\begin{Lemma}   \label{lemma2}
Let $A$ be a unital associative commutative algebra, and $\fk$ be a 
finite dimensional simple Lie algebra, both over an algebraically closed 
field $\K$. 

Any maximal ideal $J$ of the current algebra $A\otimes\fk$ must be 
of the form $I\otimes\fk$ for some maximal ideal $I$ of $A$.
\end{Lemma}

\begin{Proof}
Choose a Cartan subalgebra $\fh$ of $\fk$, and consider the Cartan 
decomposition $\fk=\fh\oplus(\oplus_{\alpha\in\Phi}\fk_{\alpha})$, where 
$\Phi$ denotes the root system of $\fk$ and as usual, 
$$\fk_{\alpha}\,=\,\{x\in\fk\,|\,\forall h\in\fh:\,\,[h,x]=\alpha(h)x\}.$$
Let $\fk$ be of rank $l$.

Let $(h_i)_{i=1}^l$ be a basis of $\fh$, and complete it into a basis 
$(x_i)_{i=1}^n$ of $\fk$ which is adapted to the Cartan decomposition, 
i.e. consisting of a basis of $\fh$ and elements from $\fk_{\alpha}$, such 
that $x_i=h_i$ for $i=1,\ldots,l$. 

We claim that for any $i$, the projection onto $x_i$, parallel (with respect 
to the Killing form) to the subspace 
generated by the other $x_j$, is the sum of compositions of inner derivations.

Indeed, for a root $\alpha$ and an element $h\in\fh$, the eigenvalues of 
$\ad(h)$ are $\beta(h)$, $\beta$ a root, or $0$. As the endomorphism $\ad(h)$ 
of $\fk$ is semi-simple, the minimal polynomial is written
$$\mu_{\ad(h)}(X)=(X-\lambda_1)\cdot\ldots\cdot(X-\lambda_r),$$
and the projectors $p^h_i$ onto the eigenspaces 
$\Ker(\ad(h)-\lambda_i)$ are 
polynomials in $\ad(h)$. More precisely, denoting 
$$Q_i(X)=(X-\lambda_1)\ldots(X-\lambda_{i-1})(X-\lambda_{i+1})\ldots
(X-\lambda_r),$$
Bezout's theorem implies that there exist $R_1,\ldots,R_r\in\C[X]$ such that 
$$R_1Q_1+\ldots R_rQ_r\,=\,1,$$
and the projector onto $\Ker(\ad(h)-\lambda_{i})$ is then $R_iQ_i(\ad(h))$.

Thus the projector $p_{\alpha}^h$ onto the eigenspace
corresponding to the eigenvalue $\alpha(h)$ is a polynomial in $\ad(h)$.
Now the composition $p_i:=p_{\alpha}^{h_1}\circ\ldots\circ p_{\alpha}^{h_l}$ is
the projection onto the $x_i$ corresponding to the root $\alpha$. 
Taking $0$ instead of $\alpha$, one obtains a polynomial in $\ad(h)$ which
represents the projector onto $\fh$. 

Now let $(\epsilon_i)_{i=1}^n$ be the
dual basis of $(x_i)_{i=1}^n$. Denote by $P_i(J)$ the image of the 
restriction of the linear maps 
${\rm id}_A\otimes(\epsilon_i\circ p_i):A\otimes\fk\to A$ to $J$. 
Let $P(J)$ be the sum of the $P_i(J)$. $P(J)$ is an ideal, because 
$J$ is closed under multiplication by an $a\in A$ in the first factor.

Indeed, the Casimir operator $\Omega=\sum_{i=1}^nx_ix^i$ where $x^i\in\fk$
for $i=1,\ldots,n$ are the elements in orthonormality to the $x_i$ with 
respect to the Killing form, satisfies $\ad\Omega=\lambda\id_{\fh}$, and 
one can renormalize the $x^i$ such that $\lambda=1$. Then one easily computes 
that
$$ab\otimes y\,=\,\sum_{i=1}^n(a\otimes\ad(x_i))(1\otimes\ad(x^i))
(b\otimes y),$$
meaning that the multiplication by $a\in A$ in the first factor can be 
expressed by a composition of inner derivations.
 
Now let us prove that $J=P(J)\otimes\fk$. 
Indeed, the inclusion $J\subset P(J)\otimes\fk$ is trivial. On the other hand, 
one has $(\id_A\otimes p_i)(J)\subset J$, because $p_i$ is given 
by a sum of compositions of inner derivations, and therefore 
$P_i(J)\otimes\fk\subset J$, which sum up to $P(J)\otimes\fk\subset J$.

This ends the proof.
\end{Proof}

\begin{Corollary} 
$$\Aut({\mathcal O}(\Sigma_k,{\mathfrak k}))\cong \Aut(\Sigma_k)
\ltimes{\mathcal O}(\Sigma_k,\Aut({\mathfrak k})).$$
\end{Corollary}

The automorphism group of the Lie algebra determines the connected component
of the identity of the automorphism group of the corresponding Lie group: 

\begin{Corollary}  \label{Aut(G)}
Let $K$ be a finite-dimensional 1-connected complex Lie group. Then
$$\Aut^0({\mathcal O}(\Sigma_k,K))\cong \Aut(\Sigma_k)
\ltimes{\mathcal O}(\Sigma_k,\Aut^0(K)),$$
where $\Aut^0(H)$ denotes the connected component of the identity in the 
automorphism group $\Aut(H)$ of the Lie group $H$.  
\end{Corollary}

\begin{Proof}
For a regular, connected and 1-connected Lie group $G$, 
every automorphism of ${\mathfrak g}$ integrates uniquely to an 
automorphism of $G$ (see \cite{Mil} Thm. 8.1). This is the case for
the regular Fr\'echet-Lie group $G={\mathcal O}(\Sigma_k,K)$ for 
finite dimensional 1-connected $K$. 
\end{Proof}

\begin{Remark} \label{Remark_Aut(G)}
\begin{enumerate}
\item[(a)] Observe that for a finite dimensional complex simple Lie algebra
${\mathfrak k}$, the group $\Aut({\mathfrak k})$ is the semi-direct product of 
the inner automorphisms ${\rm Inn}({\mathfrak k})$ of ${\mathfrak k}$ 
and the graph automorphisms (see 
\cite{Hum} Ch. 16.5, p. 87). Therefore, for $g>0$ or $k>2$, the only 
continuous families of 
automorphisms of $\Aut({\mathcal O}(\Sigma_k,{\mathfrak k}))$ are holomorphic 
maps on $\Sigma_k$ with values in ${\rm Inn}({\mathfrak k})$.
In other words, up to finite subgroups, the automorphisms of 
${\mathcal O}(\Sigma_k,{\mathfrak k})$ for $g>0$ or $k>2$ are holomorphic 
maps on $\Sigma_k$ with values in ${\rm Inn}({\mathfrak k})$.
\item[(b)] Up to finite groups, there are thus no outer automorphisms of 
${\mathcal O}(\Sigma_k,{\mathfrak k})$ for $g>0$ or $k>2$. 
This can be explained in the following 
way. The outer derivations of ${\mathcal O}(\Sigma_k,{\mathfrak k})$ are given 
by the Lie algebra $Hol(\Sigma_k)$ of holomorphic vector fields on $\Sigma_k$.
As these do not integrate into an infinite-dimensional complex Lie group, the
corresponding automorphism group is missing/finite dimensional.    
\end{enumerate}
\end{Remark}

\begin{Proposition}
If $G$ is simply connected and $\omega$ a continuous $2$-cocycle defining
a central extension $\hat{\mathfrak g}$ of ${\mathfrak g}=L(G)$ by 
${\mathfrak z}$ and $\hat{G}$ a central extension of $G$ by $Z$ integrating 
$\hat{\mathfrak g}$, then $\gamma=(\gamma_G,\gamma_Z)\in\Aut(G)\times\Aut(Z)$
lifts to an automorphism $\hat{\gamma}$ of $\hat{G}$ (fixing $Z$) if and only 
if $[\gamma\cdot\omega]=[\omega]$, i.e. if the corresponding automorphism of
${\mathfrak g}$ lifts to $\hat{\mathfrak g}$. 
\end{Proposition}

\begin{Proof}
This is Proposition V.4 in \cite{NeeMai}.
\end{Proof}

Let us investigate the condition $[\gamma\cdot\omega]=[\omega]$ for some class
of automorphisms. The central extension 
$\widehat{G}=\widehat{{\mathcal O}(\Sigma_k,K)}$ of 
$G={\mathcal O}(\Sigma_k,K)$ 
is given by the cocycle (\ref{cocycle}) in Section 
\ref{section_central_extension}, $K$ is finite dimensional 
and 1-connected with simple Lie algebra ${\mathfrak k}$ whose Killing form 
is denoted by $\kappa$. The center is 
$Z=\Big(\Omega^1(\Sigma_k)\,/\,d{\mathcal O}(\Sigma_k)\Big)\,/\,\Pi_{\Sigma_k}$. 
$\Aut(G)$ is 
described by Corollary \ref{Aut(G)}. 

\begin{Proposition}
Consider an automorphism of ${\mathcal O}(\Sigma_k,K)$ of the form 
$\gamma=(\gamma_G,\id_Z)\in\Aut(G)\times\Aut(Z)$ with 
$\gamma_G\in{\mathcal O}(\Sigma_k,\Aut(K))$.
Then $[\gamma\cdot\omega_{\Sigma_k}]=[\omega_{\Sigma_k}]$.  
\end{Proposition}

\begin{Proof}
This follows from the fact that the Killing form is invariant under 
automorphisms, because of the form of the cocycle $\omega_{\Sigma_k}$. 
\end{Proof}

\begin{Remark}
\begin{enumerate}
\item[(a)] It follows from the fact that a Lie algebra acts trivially on 
its cohomology that the identity component of the subgroup of
inner automorphisms of the Lie group always act trivially on any Lie algebra 
cocycle. This is easily shown directly using the Cartan formula for 
the action of an element $X$ in a Lie algebra ${\mathfrak g}$ on a class
$[\omega]\in H^n({\mathfrak g})$:
$$L_X\omega\,=\,(d\circ i_X+i_X\circ d)\omega\,=\,d(i_X\omega),$$
because the representative $\omega$ is closed.
\item[(b)] We believe it to be an interesting open problem to determine 
explicitly all
automorphism groups $\Aut(\widehat{{\mathcal O}(\Sigma_k,K)})$ for 
different $k\geq 1$ and different finite-dimensional 1-connected simple 
Lie groups $K$ where the extension is defined using the cocycle 
(\ref{cocycle}).
\end{enumerate}
\end{Remark} 

\section{Coadjoint Orbits}

In this section, we study the orbits of the coadjoint action of 
a central extension of $G:={\mathcal O}(\Sigma_k,K)$ on the dual space of a 
central extension of ${\mathfrak g}:={\mathcal O}(\Sigma_k,{\mathfrak k})$
with ${\mathfrak k}$ complex simple. 
As usual, we will 
restrict the dual to the so-called {\it smooth dual}. The orbits we obtain in 
this way are closely related to the flat bundles described in Section $1$. 
The main result 
is the description of the orbits in terms of conjugacy classes.

\subsection{The smooth dual}

We will define here the smooth dual of a central extension of
${\mathcal O}(\Sigma_k,{\mathfrak k})$. 
Actually, it would be more accurate to call it the holomorphic dual, but we will
stick to the classical term. Following \cite{EtiFre}, we will not consider
the universal central extension 
$\widehat{{\mathcal O}(\Sigma_k,{\mathfrak k})}$,
but only a central extension with $1$-dimensional central term. 

For this,
let us fix a singular cycle $\sigma\in Z_1(\Sigma_k)$ which we suppose to be
a simple smooth curve and which we suppose to be one of the $\ell$ generators of
$\pi_1(\Sigma_k)$. Recall the {\it compact dual}\quad 
$\eta_{\sigma}\in\Omega^1_{\rm sm,\,c}(\Sigma_k)$, where 
$\Omega^1_{\rm sm,\,c}(\Sigma_k)$ is the space of smooth differential $1$-forms
with compact support, cf. \cite{BotTu} p. 51. Here $\eta_{\sigma}$ is by 
definition a
compactly supported differential $1$-form which satisfies for any smooth 
$1$-form $\alpha$ on $\Sigma_k$
$$\int_{\sigma}\alpha\,=\,\int_{\Sigma_k}\eta_{\sigma}\wedge\alpha.$$      
We define a central extension 
$$\hat{\mathfrak g}_{\sigma}\,:=\,\widehat{{\mathcal O}(\Sigma_k,{\mathfrak k})
}_{\sigma}$$ 
of ${\mathcal O}(\Sigma_k,{\mathfrak k})$ with a 
$1$-dimensional central term by the cocycle $\Omega_{\sigma}$ given by
$$\Omega_{\sigma}(X,Y)\,:=\,\int_{\Sigma_k}\eta_{\sigma}\wedge\kappa(X,dY)\,=\,
\int_{\sigma}\kappa(X,dY),$$
where $\kappa$ is the Killing form of the simple Lie algebra ${\mathfrak k}$,
and $X,Y\in {\mathfrak g}={\mathcal O}(\Sigma_k,{\mathfrak k})$.
Elements of the central extension $\hat{\mathfrak g}_{\sigma}$ will be denoted
by $\mu k+X$, where $k$ is the central element, $\mu\in\C$ and 
$X\in{\mathcal O}(\Sigma_k,{\mathfrak k})$.

\begin{Definition}
Let $\hat{\mathfrak g}_{\sigma}^*$ be the Fr\'echet space of operators
$D=\lambda\partial+\xi$ where $\lambda\in\C$ and 
$\xi\in\Omega^1(\Sigma_k,{\mathfrak k})$ is a holomorphic $1$-form on 
$\Sigma_k$. We call $\hat{\mathfrak g}_{\sigma}^*$ the smooth dual of 
$\hat{\mathfrak g}_{\sigma}$.
\end{Definition}

In order to justify the implicit duality claim in the name of 
$\hat{\mathfrak g}_{\sigma}^*$, let us introduce the following bilinear form:

$$(\lambda\partial+\xi,\mu k +X)\,:=\,\lambda\mu + \int_{\Sigma_k}\eta_{\sigma}\wedge
\kappa(\xi,X)\,=\, \lambda\mu + \int_{\sigma}\kappa(\xi,X).$$ 

The Lie group $G={\mathcal O}(\Sigma_k,K)$ acts on 
$\hat{\mathfrak g}_{\sigma}^*$ by the following action, which is close to the
gauge action (\ref{gauge_action}) in Section \ref{section_flat_bundles}: 
$$f\cdot(\lambda\partial+\xi)\,:=\,\lambda\partial + \delta(f) + 
{\rm Ad}(f)^{-1}(\xi),$$
where $f\in G$. 
We will call this action the {\it coadjoint action} of $G$ on 
$\hat{\mathfrak g}_{\sigma}^*$. It identifies to the coadjoint action using
the duality given by the following proposition:

\begin{Proposition}
The bilinear form 
$$(-,-):\hat{\mathfrak g}_{\sigma}^*\times\hat{\mathfrak g}_{\sigma}\to\C$$
is a non-degenerate $G$-invariant pairing.
\end{Proposition}

\begin{Proof}
The form is non-degenerate: Indeed, consider for a given $\xi\not=0$, 
$\int_{\sigma}\kappa(\xi,X)$ inserting an arbitrary $X$. Then the support 
$\supp(\xi)$ has a non-empty intersection with $\im(\sigma)$. 
Working locally in a suitable open set, $\xi$ may be regarded as a 
${\mathfrak k}$-valued holomorphic function on $\Sigma_k$. Express $\xi$ as 
$\xi(z)=\sum_if_i(z)\otimes x_i$ for holomorphic functions $f_i$ and some 
elements $x_i\in{\mathfrak k}$. One may suppose that $(x_i)_{i=1}^N$ forms a basis 
of ${\mathfrak k}$ and that $x_{i_0}$ has a non-zero coefficient function in some 
open set. By non-degeneracy of the Killing form $\kappa$,
one can find $y\in{\mathfrak k}$ with $\kappa(x_{i_0},y)\not= 0$ and 
$\kappa(x_i,y)= 0$ for all $i\not= i_0$.  
Then take $X=1\otimes y$. It is clear that $(\xi,X)\not=0$.  
This reasoning works for both arguments.

The form is $G$-invariant: Indeed, the adjoint action of $G$ on 
$\hat{\mathfrak g}_{\sigma}$ is given by
$$f\cdot(\mu k+X)\,:=\,\mu k - \int_{\Sigma_k}\eta_{\sigma}\wedge\kappa(
\delta(f),X) + {\rm Ad}(f)(X).$$
With this definition, one easily computes that
$$(f\cdot(\lambda\partial+\xi),(\mu k+X))\,=\,(\lambda\partial+\xi,f\cdot
(\mu k+X)),$$
using that the Killing form is invariant under automorphisms.     
\end{Proof}

\subsection{Coadjoint orbits and flat bundles}

In this subsection, we explain the link between the framework of our Section 1 
and Section 3 in \cite{EtiFre} which establishes a bijective correspondence
between coadjoint orbits and flat bundles 
(Corollary \ref{corr_orbits_bundles}).

The coadjoint action of $G$ on $\hat{\mathfrak g}_{\sigma}^*$, introduced in 
the previous section, leaves stable the hyperplanes 
$$H_{\lambda}\,:=\,\{\,\lambda={\rm constant}\,\}\,\subset\,
\hat{\mathfrak g}_{\sigma}^*.$$
We will therefore restrict to $H_{\lambda}$ in the following. Now compare the
coadjoint action to the gauge action (\ref{gauge_action}): In $H_{\lambda}$,
these two actions coincide, thus the coadjoint orbits in 
$H_{\lambda}$ correpond to the gauge orbits in 
$Z^1_{\rm dR}(\Sigma_k,{\mathfrak k})$. The set 
$H^1_{\rm dR}(\Sigma_k,{\mathfrak k})$ of gauge orbits in 
$Z^1_{\rm dR}(\Sigma_k,{\mathfrak k})$ has been examined in Section $1$ and
identifies for $\Sigma_k$ to the space of flat $K$-principal bundles
which identifies in turn to $K^{\ell}\,/\,K$ (Corollary \ref{gauge_orbits}).

Let us show this here once again from the point of view of \cite{EtiFre},
Section $3$. 

\begin{Proposition}
Let ${\mathcal U}=\{U_i\}$ be a good open cover of $\Sigma_k$.
Let $D=\lambda\partial+\xi\in H_{\lambda}$ and consider the equation
\begin{equation}  \label{transition_functions}
\lambda(\partial\psi_i)\psi_i^{-1}+\xi\,=\,0.
\end{equation} 
Then solutions $\psi_i$ exist on $U_i$ and the functions 
$\phi_{ij}:=\psi_i^{-1}\psi_j:U_{ij}\to K$ determine transition functions of
a holomorphic $K$-principal bundle on $\Sigma_k$. Furthermore, the functions 
$\phi_{ij}$ are constant, thus the bundle is flat. 
\end{Proposition}

\begin{Proof}
Observe that the holomorphic Maurer-Cartan equation is trivially satisfied on 
$\Sigma_k$ as it is of complex dimension $1$. 
Therefore solutions to equation (\ref{transition_functions}) 
exist thanks to Theorem \ref{fundamental_theorem}.

The $\phi_{ij}$ satisfy the \v{C}ech-cocycle condition on triple intersections
$U_{ijk}:=U_i\cap U_j\cap U_k$
$$\phi_{ij}\phi_{jk}\phi_{ki}\,=\,\psi_i^{-1}\psi_j\psi_j^{-1}\psi_k\psi_k^{-1}
\psi_i\,=\,1.$$
Therefore, they determine a holomorphic $K$-principal bundle. Let us show that the 
$\phi_{ij}$ are not only holomorphic (like in \cite{EtiFre}), but also 
anti-holomorphic, and thus constant. Indeed, the transition functions
$\phi_{ij}$ are also anti-holomorphic:
\begin{eqnarray*}
\partial(\psi_i^{-1}\psi_j)&=& -\psi_i^{-1}(\partial\psi_i)\psi_i^{-1}\psi_j
+ \psi_i^{-1}(\partial\psi_j) \\
&=& \psi_i^{-1}\xi\psi_j + \psi_i^{-1}(-\xi)\psi_j\,=\,0
\end{eqnarray*}
\end{Proof}

\begin{Corollary} \label{corr_orbits_bundles}
There is a bijective correspondence between coadjoint orbits in $H_{\lambda}$
and gauge equivalence classes of flat connections on the trivial holomorphic 
$K$-principal bundle 
$K\times\Sigma_k$. 
\end{Corollary}

\begin{Proof}
Indeed, the previous proposition tells us that to a coadjoint orbits, we may 
associate an equivalence class of flat holomorphic $K$-principal bundles. 
But all these are trivial on
$\Sigma_k$ by Proposition \ref{all_trivial}. On the other hand, a flat trivial
principal $K$-bundle has a flat connection $\psi_i$ associated to the good 
cover ${\mathcal U}$, and $\partial(\psi_i^{-1}\psi_j)=0$ implies that the 
holomorphic $1$-form $\xi$ defined by $\xi=-\lambda(\partial\psi_i)\psi_i^{-1}$
is global:
\begin{eqnarray*}
\partial(\psi_i^{-1}\psi_j)&=&-\psi_i^{-1}(\partial\psi_i)\psi_i^{-1}\psi_j
+ \psi_i^{-1}(\partial\psi_j)\,=\,0 \\
&\Rightarrow&(\partial\psi_i)\psi_i^{-1}\psi_j\,=\,\partial\psi_j \\
&\Rightarrow&(\partial\psi_i)\psi_i^{-1} \,=\,(\partial\psi_j)\psi_j^{-1}.
\end{eqnarray*} 
\end{Proof}  

\begin{Remark}
For the sake of this remark, let us take $K=\C^*$. This group is not simple, 
but this will not disturb us here. We want to underline in this remark the idea
that Etingof-Frenkel's construction of the holomorphic bundle above does
not only construct a holomorphic bundle, but a holomorphic bundle {\it with a
connection} $\psi$. We choose $K=\C^*$ in order to tie this construction to
Deligne cohomology. Indeed, recall Theorem 2.2.12 of \cite{Bry}:

\begin{Theorem}
The group of isomorphism classes of pairs $(L,\nabla)$ of a line bundle $L$ on 
$M$ with connection $\nabla$ is canonically isomorphic to the Deligne 
cohomology group 
$H^1(M,\underline{\C}_M^*\stackrel{\delta}{\to}\Omega^1_{\rm sm,M,\C})$.
\end{Theorem}

As stated, the theorem is true in the smooth setting, but it can be easily
adapted to the holomorphic setting on a Stein manifold. 

Our observation is that $(\xi,\phi_{ij})$ forms a cocycle in the bicomplex of
Deligne cohomology. We believe that the part of the bicomplex which is used to
compute the Deligne $H^1$ makes also sense for a non-abelian Lie group $K$,
while the entire bicomplex makes sense only in the abelian context, i.e.
in case $K$ is $\C^*$ or more generally any abelian complex Lie group.
\end{Remark} 

\subsection{Description of coadjoint orbits using stabilizers}

In this subsection, we adapt the idea of describing orbits by stabilizers like
in \cite{MohWen} to our setting.

Let us denote by ${\mathcal O}_{(\xi,\lambda)}$ the orbit of the coadjoint 
action of $G:={\mathcal O}(\Sigma_k,K)$ on ${\mathfrak g}_{\sigma}^*$ 
passing through $\lambda\partial+\xi$. Recall from the previous subsection 
that the set of orbits in the hyperplane $H_{\lambda}$ may be identified to
$K^{\ell}\,/\,K$, where $\ell$ is the first Betti number of $\Sigma_k$ and 
$K$ acts on $K^{\ell}$ by conjugation in each factor.

\begin{Remark}  \label{Frenkel_construction}
Let us recall Frenkel's construction \cite{Fre} of the stabilizer of an
orbit of the coadjoint action for the circle $\bS^1$. For an orbit 
${\mathcal O}_{(X,\lambda)}\subset L{\mathfrak k}\oplus\C$ with 
$\lambda\not=0$, where 
$L{\mathfrak k}={\mathcal C}^{\infty}(\bS^1,{\mathfrak k})$ is the loop 
algebra, we can solve the differential equation
$$z'\,=\,-\frac{1}{\lambda}Xz.$$
Let $z_{(X,\lambda)}$ be its unique solution with initial condition 
$z_{(X,\lambda)}(0)=1\in K$. This is a path in $K$ starting at $1\in K$. 
Now, since $X\in L{\mathfrak k}$ is periodic, we get 
$z_{(X,\lambda)}(\theta+2\pi)=z_{(X,\lambda)}(\theta)z_{(X,\lambda)}(2\pi)$.
Let $\lambda\partial+Y$ be another element in the coadjoint orbit
${\mathcal O}_{(X,\lambda)}$. By the expression of the coadjoint action,
$X$ and $Y$ are linked by $Y=gXg^{-1}-\lambda\delta(g)$ for some $g$ in the 
loop group $L(K)$.
Proposition 3.2.5 in \cite{Fre} shows that the periodicity of solutions
implies the following conjugation formula
$$z_{(X,\lambda)}(\theta)\,=\,g(\theta)z_{(X,\lambda)}(\theta)g(0)^{-1}.$$ 
Now since $g$ is also periodic, $z_{(X,\lambda)}(2\pi)$ and 
$z_{(Y,\lambda)}(2\pi)$ lie in the same conjugacy class in $K$. 
Furthermore, the stabilizer of $(X,\lambda)$ in $L(K)$ is isomorphic to the 
stabilizer of $z_{(X,\lambda)}(2\pi)$ in $K$, so that we get 
$${\mathcal O}_{(X,\lambda)}\,=\,L(K)\,/\,{\rm Stab}_K(z_{(X,\lambda)}(2\pi)).
$$
An explicit description of the stabilizers is contained in \cite{MohWen}.
\end{Remark}

The idea is here to use this known description of coadjoint orbits in the loop 
group setting by restriction our holomorphic currents to circles.
For this, recall further the Lie group
homomorphism from Lemma \ref{homomorphism}
$$\alpha_K:G\to{\mathcal C}^{\infty}(\bS^1,K),$$
which is associated to $\alpha:\bS^1\to\Sigma_k$. In the following, we will 
denote ${\mathcal C}^{\infty}(\bS^1,K)$ simply by $LK=K^{\bS^1}$.
Applying the construction of $\alpha_K$ 
to the $\ell$ embedded circles in $\Sigma_k$, denoted by $\alpha_i$,
which generate $\pi_1(\Sigma_k)$, we obtain a homomorphism
$$\phi:G\to\Pi_{i=1}^{\ell}K^{\bS^1}.$$

The following proposition describes the image of the orbit 
${\mathcal O}_{(\xi,\lambda)}$, seen as $G\,/\,{\rm Stab}(\xi,\lambda)$, under 
the map induced by $\phi$:
  
\begin{Proposition}
The map $\phi$ induces an injection 
$$\bar{\phi}:G\,/\,{\rm Stab}(\xi,\lambda)\to \Pi_{i=1}^{\ell}K^{\bS^1}\,/\,
\Pi_{i=1}^{\ell} K_{C_i},$$
where $(C_1,\ldots,C_{\ell})\in\Hom(\pi_1(\Sigma_k),K)\,/\,K\cong 
K^{\ell}\,/\,K$ is the
image of $\lambda\partial+\xi$ under the period map $P$ (also called monodromy 
map in this context) and $K_{C_i}$ is the
stabilizer of the conjugacy class $C_i$ in $K$, i.e. we have 
$K\,/\,K_{C_i}=C_i$.
\end{Proposition}

\begin{Proof}
\noindent(1)\quad Definition of the map $\bar{\phi}$.

We have to show that the map $\bar{\phi}$ in the following diagram is 
well-defined.

\vspace{.5cm}
\hspace{4.5cm}
\xymatrix{
G \ar[r]^{\phi} \ar@{->>}[d] & \left(K^{\bS^1}\right)^{\ell} \ar@{->>}[d]  \\ 
G\,/\,{\rm Stab}(\xi,\lambda) \ar@{.>}[r]^{\bar{\phi}} & 
\left(K^{\bS^1}\right)^l\,/\,\Pi_{i=1}^{\ell} K_{C_i}}
\vspace{.5cm} 

We have to show that for all $g\in{\rm Stab}(\xi,\lambda)$, 
$\phi(g)\in\Pi_{i=1}^{\ell} K_{C_i}$. Now, $g\in{\rm Stab}(\xi,\lambda)$ means 
that 
$$(\lambda\partial+\xi)\cdot g\,:=\,\lambda\partial+
(\partial g)g^{-1} + {\rm Ad}\,g^{-1}(\xi)\,
\stackrel{!}{=}\,\lambda\partial+\xi.$$
In other words, 
$$\xi\,=\,(\partial g)g^{-1} + {\rm Ad}\,g^{-1}(\xi)\,=:\,\xi\cdot g.$$

By Remark \ref{Frenkel_construction}, we have a commutative diagram

\vspace{.5cm}
\hspace{4.5cm}
\xymatrix{
H^1_{\rm dR}(\Sigma_k,{\mathfrak k}) \ar[r]^{\hspace{-1cm}P} 
\ar[d]^{E:=|_{\alpha_i}} & 
H^1(\pi_1(\Sigma_k),K)=K^{\ell}\,/\,K \ar@{->>}[d]^{\pi_i}  \\ 
H^1_{\rm dR, sm}(\bS^1,{\mathfrak k}) \ar[r]^{F} & 
K\,/\,K}
\vspace{.5cm} 

Here the map $|_{\alpha_i}:H^1_{\rm dR}(\Sigma_k,{\mathfrak k})\to
H^1_{\rm dR, sm}(\bS^1,{\mathfrak k})$ restricts gauge orbits of forms on
$\Sigma_k$ to gauge orbits of forms on the $i$th circle 
$\alpha_i\subset\Sigma_k$. The map $P$ is the period map from Section $1$ and
$\pi_i$ is the
projection onto the $i$th factor in the product. The map $F$ is Frenkel's 
construction and $K\,/\,K$ is the set of conjugacy classes of $K$.  

By commutativity, we have on the one hand
$$\pi_i\circ P([\xi])\,=\, C_i\,=\,F\circ E([\xi]),$$
and on the other hand
$$\pi_i\circ P([\xi\cdot g])\,=\, C_i\,=\,F\circ E([\xi\cdot g]).$$
The equality $\xi=\xi\cdot g$ implies therefore that $\xi|_{\alpha_i}$ and
$(\xi\cdot g)|_{\alpha_i}$ are sent to $C_i$ under Frenkel's construction.
We get 
$$\pi_i(\phi(g))\,=\,g|_{\alpha_i}\in{\rm Stab}(\xi|_{\alpha_i},\lambda)
\,=\,G_{C_i}.$$

\noindent(2)\quad Injectivity of the map $\bar{\phi}$.  
 
Here we have to show for $f,g\in G$ with $\bar{\phi}(f)=\bar{\phi}(g)$ that 
$fg^{-1}\in{\rm Stab}(\xi,\lambda)$. This means explicitly
$$\delta(f)+{\rm Ad}\,f^{-1}\xi\,=\,\delta(g)+{\rm Ad}\,g^{-1}\xi.$$
This equation between holomorphic functions may be shown by restriction to 
{\it one} circle $\alpha_i$, because $f,g$ are holomorphic. 
On one circle $\alpha_i$, this is exactly the outcome of Remark 
\ref{Frenkel_construction}.      
\end{Proof}

\begin{rem}
In fact, it suffices to take only one factor in the product 
$\Pi_{i=1}^{\ell}K^{\bS^1}$, but in order to treat the circles $\alpha_i$
on an equal basis, we chose to take them all. 
\end{rem}

\subsection{The symplectic form on the orbits}

It is known since work of Kostant \cite{Kos}, Kirillov \cite{Kir} and 
Souriau \cite{Sou} that coadjoint 
orbits are symplectic manifolds. Let us recall in this section the corresponding
symplectic form. The main outcome of this section is that in our setting, the 
orbits are {\it holomorphic symplectic manifolds}. 

Consider here one embedded circle $\alpha:\bS^1\to\Sigma_k$ in $\Sigma_k$.
Let us distinguish $\alpha$ from its image $\alpha(\bS^1)=:C\subset\Sigma_k$. 
Denote by 
$$\phi:{\mathcal O}(\Sigma_k,K)\to{\mathcal C}^{\infty}(\bS^1,K)$$
the restriction map $f\mapsto f\circ\alpha$ to the circle $C$. This is 
a homomorphism of complex Fr\'echet-Lie groups. We have seen above that 
$\phi$ induces a map 
$$\bar{\phi}:{\mathcal O}_{(\xi,\lambda)}\to{\mathcal O}_{(\xi|_C,\lambda)},$$
the only new issue being that there is only one circle.   

Let us describe tangent vectors $v$ to the orbit ${\mathcal O}_{(\xi,\lambda)}$.
The vector $v\in T_{(\lambda',\xi')}{\mathcal O}_{(\xi,\lambda)}$ is described by
 a smooth curve $\gamma:[-1,1]\to{\mathcal O}_{(\xi,\lambda)}$ such that 
$\gamma(0)=(\lambda',\xi')$ and $v=\dot{\gamma}(0)$. As the orbit  
${\mathcal O}_{(\xi,\lambda)}$ lies entirely in the hyperplane 
$$H_{\lambda}\,=\,\{\,\lambda={\rm constant}\,\}\,\subset\,
\hat{\mathfrak g}_{\sigma}^*,$$
we have $\lambda'=\lambda$. The curve can be taken explicitly as
$$\gamma_v(t)\,:=\,(\lambda,\xi')\cdot e^{t\tilde{X}}$$
for some $\tilde{X}\in{\mathcal O}(\Sigma_k,{\mathfrak k})$.   

\begin{Definition}
The KKS-form $\omega_h$ on the orbit 
${\mathcal O}_{(\lambda,\xi)}$ is given for all tangent vectors 
$v,w\in T_{(\lambda,\xi')}{\mathcal O}_{(\xi,\lambda)}$ by 
$$\omega_h(v,w)\,:=\,\int_C\kappa(\xi',[\tilde{X},\tilde{Y}]),$$
where $v,w$ are represented by $(\xi',\lambda)\cdot e^{t\tilde{X}}$ resp.
$(\xi',\lambda)\cdot e^{t\tilde{Y}}$ as explained above.
\end{Definition}

As always in infinite dimensions, $\omega_h$ is a ``symplectic'' form
in the sense that it is {\it weakly symplectic}, i.e. instead of demanding 
a {\it non-degenerate} skewsymmetric bilinear form, one only demands that 
the form induces an injection into the dual space. 
The form is nonetheless supposed to be closed. These properties are shown 
below.  

In the same way, there is the standard symplectic form $\omega_s$ on the
coadjoint orbits of the loop group ${\mathcal C}^{\infty}(\bS^1,K)$. 
The relation between the two symplectic forms is simply:

\begin{Proposition}
For all tangent vectors 
$v,w\in T_{(\lambda,\xi')}{\mathcal O}_{(\xi,\lambda)}$, we have 
$$\omega_h(v,w)\,=\,\omega_s(T\bar{\phi}(v),T\bar{\phi}(w)).$$
\end{Proposition}

\begin{Proof}
We have 
\begin{eqnarray*} 
\omega_h(v,w)&=&\int_C\kappa(\xi',[\tilde{X},\tilde{Y}]) \\
&=&\int_{\bS^1}\kappa(\xi'|_C,T_{(\lambda,\xi')}\bar{\phi}
[\tilde{X},\tilde{Y}])   \\
&=&\int_{\bS^1}\kappa(\xi'|_C,[T_{(\lambda,\xi')}\bar{\phi}(\tilde{X}),
T_{(\lambda,\xi')}\bar{\phi}(\tilde{Y})])   \\
&=&\omega_s(T\bar{\phi}(v),T\bar{\phi}(w)). 
\end{eqnarray*}
\end{Proof}

\begin{Corollary}
The form $\omega_h$ is closed and weakly symplectic. 
\end{Corollary}

\begin{Proof}
The first of these two properties is inherited from $\omega_s$ thanks to 
the preceeding proposition. The second assertion follows form the 
non-degeneracy and invariance of $\kappa$, because
$$\omega_h(v,-)\,=\,\int_C\kappa(\xi',[\tilde{X},-])\,=\,
\int_C\kappa([\xi',\tilde{X}],-),$$
and this latter expression, seen as a map 
$$T_{(\lambda,\xi')}{\mathcal O}_{(\xi,\lambda)}\to(T_{(\lambda,\xi')}
{\mathcal O}_{(\xi,\lambda)})^*,\,\,\,\,\,v\mapsto\omega_h(v,-)$$
is clearly injective.
\end{Proof}

\end{document}